\documentclass[reqno,10pt,centertags,draft]{amsart}
\usepackage{amsmath,amsthm,amscd,amssymb,latexsym,upref}

\newcommand{\bbN}{{\mathbb{N}}}
\newcommand{\bbR}{{\mathbb{R}}}

\newcommand{\bbC}{{\mathbb{C}}}

\newcommand{\calR}{{\mathcal R}}

\newcommand{\dott}{\,\cdot\,}
\newcommand{\no}{\nonumber}
\newcommand{\lb}{\label}
\newcommand{\f}{\frac}

\newcommand{\ti}{\tilde  }
\newcommand{\wti}{\widetilde  }

\newcommand{\loc}{\text{\rm{loc}}}

\newcommand{\rank}{\text{\rm{rank}}}

\newcommand{\dom}{\text{\rm{dom}}}

\newcommand{\supp}{\text{\rm{supp}}}
\newcommand{\AC}{\text{\rm{AC}}}
\newcommand{\bi}{\bibitem}
\newcommand{\hatt}{\widehat}
\newcommand{\beq}{\begin{equation}}
\newcommand{\eeq}{\end{equation}}
\newcommand{\ba}{\begin{align}}
\newcommand{\ea}{\end{align}}
\newcommand{\vT}{\vartheta}
\newcommand{\vP}{\varphi}

\newcommand{\vE}{\eta}

\renewcommand{\Re}{\text{\rm Re}}
\renewcommand{\Im}{\text{\rm Im}}
\renewcommand{\ln}{\text{\rm ln}}

\allowdisplaybreaks
\numberwithin{equation}{section}

\newtheorem{theorem}{Theorem}[section]
\newtheorem{lemma}[theorem]{Lemma}
\newtheorem{corollary}[theorem]{Corollary}
\newtheorem{hypothesis}[theorem]{Hypothesis}
\theoremstyle{definition}
\newtheorem{definition}[theorem]{Definition}

\theoremstyle{remark}
\newtheorem{remark}[theorem]{Remark}

\newcommand{\abs}[1]{\lvert#1\rvert}

\begin{document}

\title[Weyl-Titchmarsh Matrix]{Weyl-Titchmarsh
$M$-Function Asymptotics  for Matrix-Valued Schr\"odinger
Operators}
\author[Clark and Gesztesy]{Steve Clark and Fritz Gesztesy}
\address{Department of Mathematics and Statistics, University
of Missouri-Rolla, Rolla, MO 65409, USA}
\email{sclark@umr.edu}
\urladdr{http://www.umr.edu/\~{ }clark}
\address{Department of Mathematics,
University of
Missouri, Columbia, MO
65211, USA}
\email{fritz@math.missouri.edu\newline
\indent{\it URL:}
http://www.math.missouri.edu/people/fgesztesy.html}
\thanks{{\it Proc. London Math. Soc.} {\bf 82}, 701--724 (2001).}

\begin{abstract}
We explicitly determine the high-energy asymptotics for
Weyl-Titchmarsh matrices associated with general
matrix-valued Schr\"odinger operators on a half-line. 
\end{abstract}

\maketitle

\section{Introduction}\lb{s1}

The high-energy asymptotics $z\to\infty$ of scalar-valued
Weyl-Titchmarsh functions $m_+(z,x_0)$ associated with general
half-line Schr\"odinger differential expressions of the type
$-\f{d^2}{dx^2}+q(x)$, with real-valued coefficients $q^{(n)}\in
L^1([x_0,c])$ for some $n\in\bbN_0(=\bbN\cup\{0\})$ and
all $c>x_0$, received enormous attention over the past three
decades as can be inferred, for instance from \cite{At81},
\cite{At88}, \cite{Be88}, \cite{Be89}, \cite{BM97}, \cite{DL91},
\cite{Ev72}, \cite{EH78}, \cite{EHS83}, \cite{Ha83}--\cite{Ha87},
\cite{Hi69}, \cite{HKS89}, \cite{KK86}, \cite{KK87}, \cite{Si98}
and the literature therein. Hence it may perhaps come as a surprise
that the corresponding matrix extension of this problem,
considering general matrix-valued differential expressions of the
type $-I_m\f{d^2}{dx^2}+Q(x)$, with $I_m$ the identity matrix in
$\bbC^m$, $m\in\bbN$, and $Q(x)$ a self-adjoint matrix satisfying
$Q^{(n)}\in L^1([x_0,c])^{m\times m}$ for some $n\in\bbN_0$ and
all $c>x_0$, received hardly any attention at all. In fact, we are
only aware of a total of two papers devoted to Weyl-Titchmarsh
$M$-function asymptotics for general half-line matrix
Schr\"odinger operators: one by Clark \cite{Cl92} in the case
$m=2$, $n=0$, and a second one, an apparently unpublished manuscript
by Atkinson \cite{At88} in the case $m\in\bbN$, $n=0$. Both authors
are solely concerned with determining the leading term in the
asymptotic high-energy expansion of the Weyl-Titchmarsh matrix
$M_+(z,x_0)$ as
$z\to\infty$ without investigating higher-order expansion
coefficients. (It should be noted that this observation discounts
papers in the special scattering theoretic case concerned with
short-range coefficients $Q^{(n)}\in L^1([x_0,\infty);
(1+|x|)dx)^{m\times m}$, where straightforward iterations of
Volterra-type integral equations yield the asymptotic high-energy
expansion of $M_+(z,x_0)$ to any order, cf.~Lemma~\ref{4.1}.) A
quick look at Atkinson's preprint reveals the nontrivial
nature of determining even the leading-order term of the
asymptotic expansion of $M_+(z,x_0)$ as $z\to\infty$,
\begin{equation}
M_+(z,x_0)\underset{z\to \infty}{=} iz^{1/2}I_m + o(|z|^{1/2})
\lb{1.1}
\end{equation}
for $z$ in sectors in the open upper half-plane due to a variety of
 variable transformations, which finally reduce the
problem to the study of a managable Riccati system of differential
equations.

Our principal motivation in studying this problem stems from our
general interest in operator-valued Herglotz functions
(cf.~\cite{Ca76}, \cite{GKMT98}, \cite{GM99}, \cite{GM99a},
\cite{GMN98}, \cite{GMT98}, \cite{GT97},
\cite{Sh71}) and its possible applications in the areas of inverse
spectral theory and completely integrable systems. More precisely,
using higher-order expansions of the type
\eqref{1.1}, one can prove trace formulas for $Q(x)$ and certain
higher-order differential polynomials in
$Q(x)$, following an approach pioneered in \cite{GS96} (see also
\cite{GH97}, \cite{GHSZ95}). These trace formulas, in
turn, then can be used to prove various results in inverse spectral
theory for matrix-valued Schr\"odinger operators
$H=-I_m\f{d^2}{dx^2}+Q$ in $L^2(\bbR)^m$. For instance, using the
principal result of this paper, Theorem~\ref{t4.6}, and its
straightforward application to the asymptotic high-energy
expansion of the diagonal Green's matrix
$G(z,x,x)=(H-z)^{-1}(x,x)$ of $H$, the following matrix-valued
extension of a classical uniqueness result of Borg \cite{Bo46} has
been obtained in \cite{CGHL99}.

\begin{theorem} [\cite{CGHL99}] \lb{t1.1} Suppose $H$ is
reflectionless {\rm (}e.g., $Q$ is periodic and H has uniform spectral
multiplicity $2m${\rm )} and spectrum $[E_0,\infty)$ for some
$E_0\in\bbR$. Then
\begin{equation}
Q(x)=E_0I_m \text{ for a.e.~$x\in\bbR$}. \lb{1.2}
\end{equation}
\end{theorem}

For related results see, for instance,
\cite{AK92},
\cite{Ca98}, \cite{Ca98a}, \cite{Ca99}, \cite{Ch99}, \cite{Ch99a},
\cite{CS97}, \cite{CGR99}, \cite{De95}, \cite{JL97}, \cite{Kr83},
\cite{Ma94}, \cite{Ma98}, \cite{Pa95}, \cite{Sa94}. Incidentally,
the higher-order differential polynomials in $Q(x)$ just alluded to
represent the Korteweg-deVries (KdV) invariants (i.e.,  densities
associated with KdV conservation laws) and hence open the link to
infinite-dimensional completely integrable systems (cf.
\cite{AG98},
\cite{Ch96}, \cite{Di91}, \cite{DK97}, \cite{Du77}, \cite{Du83},
\cite{GD77}, \cite{GKS97}, \cite{Ma78}, \cite{Ma88}, \cite{MO82},
\cite{OMG81}, \cite{Sa94a}, \cite{Sc83} and the references therein).

In Section~\ref{s2} we briefly recall some of the basic notions
of Weyl-Titchmarsh theory for singular matrix Schr\"odinger
operators on a half-line as developed in detail by Hinton and
Shaw \cite{HS81}--\cite{HS86} (see also \cite{At64}, \cite{Hi69},
\cite{Jo87}, \cite{KR74}, \cite{KS88}, \cite{Kr89a}, \cite{Kr89b},
\cite{Or76},
\cite{Re96}, \cite{Ro60}, \cite{Ro69}, \cite{Sa94a},
\cite{We87}). In fact, most of these references deal with more
general singular Hamiltonian systems and hence we here specialize
this material to the matrix-valued Schr\"odinger operator case at
hand. The subsequent section is devoted to a proof of the
leading-order asymptotic high-energy expansion
\eqref{1.1} of $M_+(z,x_0)$, originally due to Atkinson. Since
the result is a fundamental ingredient for our principal
Section~\ref{s4}, and to the best of our knowledge, his manuscript
\cite{At88} remained unpublished, we provide a detailed
description of his approach in Section~\ref{s3}. Finally,
Section~\ref{s4} develops a systematic higher-order high-energy
asymptotic expansion of $M_+(z,x_0)$, combining Atkinson's result
in Section~\ref{s3} with matrix-valued extensions of some methods
based on an associated Riccati equation. More precisely,
following a technique in \cite{GS98} in the scalar-valued context,
we show how to derive the general high-energy asymptotic expansion
of
$M_+(z,x_0)$ as $z\to\infty$ by combining Atkinson's
leading-order term in \eqref{1.1} and the corresponding asymptotic
expansion of $M_+(z,x_0)$ in the special case where $Q$ has
compact support.

Analogous results for Dirac-type operators are in preparation
\cite{CG00}.

\section{Weyl-Titchmarsh Matrices}
\lb{s2}
In this section we briefly recall the Weyl--Titchmarsh
theory for matrix-valued Schr\"{o}dinger
operators on half-lines $(x_0,\infty)$ for some
$x_0\in\bbR$.

Throughout this note all $m\times m$ matrices
$M\in\bbC^{m\times m}$ will be  considered over
the field of complex numbers $\bbC$. Moreover, $I_p$ denotes
the identity matrix in $\bbC^{p\times p}$ for $p\in\bbN$, $M^*$
the adjoint (i.e., complex
conjugate transpose), $M^t$ the transpose of the matrix
$M$, and $\AC([a,b])$ denotes the
set of absolutely continuous functions on $[a,b]$.

The basic assumption of this paper will be the following.

\begin{hypothesis}\lb{h2.1}
Let $m\in\bbN$ and suppose
$Q=Q^*\in L^1([x_0,c])^{m\times m}$ for all $c>x_0$.
\end{hypothesis}

Given Hypothesis~\ref{h2.1} we consider the matrix-valued
Schr\"odinger operator $H_+$ in $L^2((x_0,\infty))^m$,
\begin{align}
&H_+=-I_m\f{d^2}{dx^2} + Q, \lb{2.1} \\
&\dom(H_+)=\{g\in L^2((x_0,\infty))^m \,|\,g,g'\in
AC([x_0,c])^m
\text{ for all } c>x_0; \lim_{x\downarrow x_0} g(x)=0, \no \\
& \hspace*{5.4cm} \text{s.-a.b.c. at }\infty; (-I_mg''+Qg)\in
L^2((x_0,\infty))^m \}, \no
\end{align}
where the abbreviation s.-a.b.c. denotes a fixed
self-adjoint boundary condition at $\infty$ throughout this
paper if the differential expression
$-I_m\f{d^2}{dx^2} + Q(x)$  is not in the limit point case at
$\infty$ (cf.~the paragraph following Definition~\ref{d2.7}) and is
disregarded in the case where
$-I_m\f{d^2}{dx^2} + Q(x)$ is in the limit point case at $\infty$.
Naturally  associated with the operator $H_+$ is the matrix
Schr\"odinger  equation
\begin{subequations}\lb{2.2}
\begin{equation}\lb{2.2a}
-\psi''(z,x) + Q(x)\psi(z,x) = z \psi(z,x), \quad
z\in\bbC, \, x\geq x_0,
\end{equation}
where $z$ represents a spectral parameter and $\psi(z,x)$
is assumed to satisfy
\begin{equation}\lb{2.2b}
\psi'(z,\cdot)\in AC([x_0,c])^{m\times r} \text{ for all } c>0,
\end{equation}
with the value of $r$ depending upon the context involved.
\end{subequations}
Equation \eqref{2.2a} may also be expressed as the system
\begin{subequations}\lb{2.3}
\begin{equation}\lb{2.3a}
J\varPsi (z,x)'=(zA(x)+B(x))\varPsi (z,x), \quad z\in\bbC, \, x\geq x_0,
\end{equation}
where
\begin{equation}\lb{2.3b}
\varPsi(z,x) = \begin{bmatrix}\psi(z,x)\\
\psi'(z,x) \end{bmatrix},\hspace{10pt}J=\begin{bmatrix}
0& -I_m\\ I_m
& 0 \end{bmatrix},
\end{equation}
\begin{equation}\lb{2.3c}
A(x) =\begin{bmatrix}I_m& 0\\0& 0\end{bmatrix},
\hspace{10pt}  B(x)=\begin{bmatrix}-Q(x)& 0\\0& I_m\end{bmatrix},
\end{equation}
and where
\begin{equation}\lb{2.3d}
\varPsi(z,\cdot)\in AC([x_0,c])^{2m\times r} \text{ for all } c>0;
\end{equation}
again, with the value of $r$ depending upon the context  involved.
\end{subequations}

Next we turn to Weyl-Titchmarsh theory associated
with \eqref{2.1} and briefly recall some
of the results developed by Hinton and Shaw in a series
of papers devoted to spectral theory of
(singular) Hamiltonian systems \cite{HS81}--\cite{HS86} (see
also \cite{At64}, \cite{Hi69},
\cite{Jo87}, \cite{KR74}, \cite{KS88}, \cite{Kr89a}, \cite{Kr89b},
\cite{Or76},
\cite{Re96}, \cite{Ro60}, \cite{Ro69}, \cite{Sa94a},
\cite{We87}). While they discuss
 much more general systems, we here confine ourselves to
the special case of  matrix-valued Schr\"{o}dinger
operators
governed by Hypothesis~\ref{h2.1}.

Let $\Psi(z,x,x_0)$ be a normalized fundamental matrix of
solutions of
\eqref{2.3} at some
$x_0\in\bbR$; that is, $\Psi(z,x,x_0)$ is of the
type
\begin{subequations}\lb{2.4}
\begin{align}
\Psi(z,x,x_0)&=[\psi_{j,k}(z,x,x_0)]_{j,k=1}^2 \notag\\
&= \begin{bmatrix}\Theta(z,x,x_0)& \Phi(z,x,x_0) \end{bmatrix} =
\begin{bmatrix}\theta(z,x,x_0) & \phi(z,x,x_0)\\
\theta'(z,x,x_0)& \phi'(z,x,x_0)\end{bmatrix}, \lb{2.4a} \\
\Psi(z,x_0,x_0)&=I_{2m}, \lb{2.4b}
\end{align}
where $\Theta(z,x,x_0),\ \Phi(z,x,x_0)\in \bbC^{2m\times m}$,
and where  $\theta(z,x,x_0)$ and $\phi(z,x,x_0)$ represent a
fundamental system of $m\times m$
matrix-valued solutions of \eqref{2.2}, entire with
respect to $z\in\bbC$, and normalized according to
\eqref{2.4b}; that is,
\begin{equation}
\theta(z,x_0,x_0)=\phi'(z,x_0,x_0)=I_m, \quad
\theta'(z,x_0,x_0)=\phi(z,x_0,x_0)=0. \lb{2.4c}
\end{equation}
\end{subequations}
Let $\beta=[\beta_1\ \beta_2]$; with  $\beta_1$, $\beta_2\in\bbC^{m\times
m}$.  We assume
\begin{subequations}\lb{2.5}
\begin{equation}\lb{2.5a}
\ker{(\beta_1^*)}\cap\ker{(\beta_2^*)}=\{0\},
\text{ or equivalently that  } \rank(\beta) = m,
\end{equation}
and that either
\begin{equation}\lb{2.5b}
\pm(1/2i)\beta J \beta^*=\pm\Im{(\beta_2\beta_1^*)> 0}, \quad\text{or}\quad
\Im{(\beta_2\beta_1^*) = 0}.
\end{equation}
\end{subequations}
Here we denote, as usual,
$\Im(M)=(M-M^*)/(2i)$ and $\Re(M)=(M+M^*)/2$.

We now prove the following result concerning  $\Phi(z,x,x_0)$.
\begin{lemma} \lb{l2.2}
Let   $\Phi(z,x,x_0)$ be defined in \eqref{2.4},
and let $\beta\in\bbC^{m\times 2m}$ satisfy
\eqref{2.5}. If $\beta\Phi(z,c,x_0)=\beta_1\phi(z,c,x_0)
+\beta_2\phi'(z,c,x_0)$ is singular for $c>x_0$ and
$\Im{(\beta_2\beta_1^*)}\gtrless 0$,
then $\Im{(z)}\lessgtr 0$.  If $\beta\Phi(z,c,x_0)$ is singular  for
$c>x_0$ and $\Im{(\beta_2\beta_1^*)}= 0$, then $\Im{(z)}= 0$.
\end{lemma}
\begin{proof}
Let $\varPsi(z,x)$ satisfy \eqref{2.3}. Then,
\begin{subequations}\lb{2.6}
\begin{equation}\lb{2.6a}
\varPsi(z,c)^*J\varPsi(z,c) - \varPsi(z,x_0)^*J \varPsi(z,x_0)
=2i\Im{(z)}\int_{x_0}^c ds\, \varPsi(z,s)^*A(s)\varPsi(z,s).
\end{equation}
Let $\varPsi(z,x)=\Phi(z,x,x_0)$ and note that
\begin{equation}\lb{2.6b}
\varPsi(z,x_0)^*J\varPsi(z,x_0)=[0 \ I_m]J[0 \ I_m]^t=0.
\end{equation}
Moreover, note that
\begin{equation}\lb{2.6c}
\varPsi(z,c)^*J\varPsi(z,c)=
\phi'(z,c,x_0)^*\phi(z,c,x_0)-\phi(z,c,x_0)^*\phi'(z,c,x_0).
\end{equation}
\end{subequations}

Now suppose that $\beta\Phi(z,c,x_0)v =0$ for
$v\in\bbC^m$, $v\ne 0$.  If $\Im{(\beta_2\beta_1^*)}\gtrless 0$, then
$\phi'(z,c,x_0)v=-\beta_2^{-1}\beta_1\phi(z,c,x_0)v$. Thus by
\eqref{2.6}  we obtain,
\begin{align}
&v^*\phi(z,c,x_0)^*\beta_2^{-1}[\Im{(\beta_2\beta_1^*)}]
{\beta_2^*}^{-1}
\phi(z,c,x_0)v \no \\
&= -\Im{(z)} \int_{x_0}^c ds \ v^*
\phi(z,s,x_0)^*\phi(z,s,x_0)v,
\end{align}
by which we obtain the first part of the lemma.  However, if
$ 0=\Im{(\beta_2\beta_1^*)}=(1/2i)\beta J \beta^*$ then by \eqref{2.5a}
there is a
$w\in \bbC^m$ such that
\begin{equation}
 \Phi(z,c,x_0)v =J\beta^*w.
\end{equation}
For this case, by \eqref{2.6} we obtain
\begin{equation}
0=\Im{(z)} \int_{x_0}^c ds \ v^*
\phi(z,s,x_0)^*\phi(z,s,x_0)v,
\end{equation}
which implies that $\Im{(z)}=0$.
\end{proof}
\noindent One can also prove the following result.
\begin{lemma}\lb{l2.3}
Let  $\Theta(z,x,x_0) $ and $\Phi(z,x,x_0)$ be defined in \eqref{2.4},
and let $\beta\in\bbC^{m\times 2m}$ satisfy
\eqref{2.5}. For $c>x_0$,
$\beta\Phi(z,c,x_0)$ is singular  if and only if $z$ is an
eigenvalue for the  regular  boundary value
problem given
by \eqref{2.3} together with the boundary conditions
\begin{equation}
[I_m \ 0]\varPsi(z,x_0)=0, \quad \beta\varPsi(z,c)=0.
\end{equation}
\end{lemma}
Lemmas~\ref{l2.2} and \ref{l2.3} show that with appropriate values of
$z\in\bbC$ one may define a certain $m\times m$  matrix.
\begin{definition}\lb{d2.4}
Let  $\Theta(z,x) $ and $\Phi(z,x)$ be defined in \eqref{2.4}, and let
$\beta\in\bbC^{m\times 2m}$ satisfy  \eqref{2.5}. Then, for $c>x_0$  and
for $\beta\Phi(z,c,x_0)$ nonsingular, let
\begin{equation}\lb{2.11}
M(z,c,x_0,\beta) =
-(\beta\Phi(z,c,x_0))^{-1}(\beta\Theta(z,c,x_0)).
\end{equation}
$M(z,c,x_0,\beta)$ is said to be the {\it Weyl-Titchmarsh M-function} for
the regular boundary value problem described in Lemma~\ref{l2.3}.
Also, let $D(z,c,x_0)$  denote the set of all  $M(z,c,x_0,\beta)$ defined
in \eqref{2.11} together with the assumption \eqref{2.5b}.
\end{definition}
About $M(z,c,x_0,\beta)$, we prove the following result.
\begin{lemma} \lb{l2.5}
Let $c>x_0$.  If $\Im{(\beta_2\beta_1^*)}>0$, then
$M(z,c,x_0,\beta) $ is well-defined if  $\Im{(z)}\ge 0$; if
$\Im{(\beta_2\beta_1^*)}<0$, then $M(z,c,x_0,\beta) $ is well-defined for
$\Im{(z)}\le 0$; and if $\Im{(\beta_2\beta_1^*)}=0$ then
$M(z,c,x_0,\beta) $ is well-defined for $z\in \bbC\backslash\bbR$.
Moreover, if  $\Im{(\beta_2 \beta_1^*)}\gtreqless 0 $ and
$ \Im{(z)}\gtrless 0 $, then
\begin{equation}\lb{2.12}
\Im{ (\pm M(z,c,x_0,\beta))} > 0.
\end{equation}
\end{lemma}
\begin{proof}
Except for the conditions given in \eqref{2.12}, the claims
made in the statement of this lemma are an immediate consequence of
Lemma~\ref{l2.2}.

To verify the conditions given in
\eqref{2.12},  we begin by letting
\begin{equation}\lb{2.13}
U(z,x,x_0)=
\begin{bmatrix}\Theta(z,x,x_0)& \Phi(z,x,x_0) \end{bmatrix}\begin{bmatrix}
I_m\\
M \end{bmatrix},
\end{equation}
where $\Theta(z,x,x_0) $ and $\Phi(z,x,x_0)$ are
defined in \eqref{2.4}, and where $M\in\bbC^{m\times m}$.
Now,  let $\Psi(z,x)=U(z,x,x_0) $  in \eqref{2.6a} and note that
\begin{equation}\lb{2.14}
\Psi(z,x_0)^*J\Psi(z,x_0)= 2i\Im (M).
\end{equation}
Moreover, if $M=M(z,c,x_0,\beta)$,  then also note that
\begin{equation}\lb{2.15}
\beta \Psi(z,c)=\beta U(z,c,x_0)=0.
\end{equation}

If $(1/2i)\beta J \beta^*=\Im{(\beta_2\beta_1^*)}=0$, then by \eqref{2.5a}
and \eqref{2.15}
there is a matrix $C\in\bbC^{m\times m}$ such that $\Psi(z,c)=
U(z,c,x_0)=J\beta^*C$. Thus,
\begin{equation} \lb{2.16}
\Psi(z,c)^*J\Psi(z,c)= C^* \beta J \beta^*C =0,
\end{equation}
and  \eqref{2.12} follows immediately from \eqref{2.14}, \eqref{2.16} and
\eqref{2.6a} when $\Im{(z)}\ne 0$.
On the other hand, if $\Im{(\beta_2\beta_1^*)}\gtrless 0$,  then
\eqref{2.15}   implies that $u'(z,c) = -\beta_2^{-1}\beta_1 u(z,c)$. As a
result,
\begin{equation}\lb{2.17}
 \Psi(z,c)^*J \Psi(z,c) =
-2iu(z,c)^*\beta_2^{-1}\Im{(\beta_2\beta_1^*)}{\beta_2^*}^{-1}u(z,c),
\end{equation}
and again \eqref{2.12} follows from \eqref{2.14}, \eqref{2.17}, and
\eqref{2.6a} when $\Im{(z)}\ne 0$.
\end{proof}
\begin{remark}\lb{r2.6}
By Lemma~\ref{l2.5}, $M(z,c,x_0,\beta)$ is a Herglotz
function of rank $m$. Moreover, $D(z,c,x_0)$ is 
contained in the standard {\em Weyl disk} of matrices associated  
with \eqref{2.2} or
\eqref{2.3} for fixed values of $z$,
$c$, and $x_0$ (cf.~\cite{CG00}).
These disks nest with respect to increasing
values of $c$; that is, $D(z,c_2,x_0)\subseteq D(z,c_1,x_0)$
whenever $x_0<c_1\le c_2$.  This leads to the existence of a closed
 limiting set
$D_+(z,x_0)$, which we call a limiting disk and which
corresponds to the standard limiting Weyl disk 
(cf.~\cite{CG00}, \cite{HS81}, \cite{HS83},
\cite{HS84}, \cite{Or76}).
\end{remark}

\begin{definition}\lb{d2.7}
Let $ D_{+}(z,x_0) $ denote the limit, as $c\to \infty$, of the nested
collection of sets $\overline{D(z,c,x_0)}$ given in Definition~\ref{d2.4},
and let $M_+(z,x_0)$ denote a point 
of the limit disk $D_+(z,x_0)$
associated with \eqref{2.2}.
\end{definition}

We note that  $D_+(z,x_0)$ consists of a unique
matrix  if and only if
$-I_md^2/dx^2+Q(x)$ is in the limit point case at $+\infty$ and that
$D_+(z,x_0)$ consists of a closed, convex, and bounded set with 
nonempty interior if
and only if  $-I_md^2/dx^2+Q(x)$ is in the limit circle case at $+\infty$.
In both of these cases, we note that elements of  $\partial D_+(z,x_0)$
represent {\em half-line Weyl-Titchmarsh  matrices} where  each such
element is associated with the construction of self-adjoint operator $H_+$
on the half-line.  However, for those operators  that fall between the
limit point and limit circle cases, Hinton and Schneider have noted that
not every element of $\partial D_+(z,x_0)$ is a half-line Weyl-Titchmarsh
matrix and have characterized those elements of the boundary that are
(cf.~\cite{HSH93}, \cite{HSH97}).

For convenience of the reader we quickly summarize the principal
results on half-line Weyl-Titchmarsh  matrices in the following
theorem.
\begin{theorem}[\cite{AD56}, \cite{Ca76}, \cite{GT97},
\cite{HS81}, \cite{HS82}, \cite{HS86},
\cite{KS88}] \lb{t2.8}
Assume Hypothesis~\ref{h2.1}, $x_0\in\bbR$,
 $z\in\bbC\backslash\bbR$, and let $M_+(z,x_0)$ be a half-line
Weyl-Titchmarsh matrix. Then \\
(i) $M_+(z,x_0)$ is a matrix-valued Herglotz
function of maximal rank. In particular,
\begin{gather}
\Im(M_+(z,x_0)) > 0, \quad z\in\bbC_+,
\lb{2.18} \\
M_+(\bar z,x_0)=M_+( z,x_0)^*,\lb{2.19} \\
\rank (M_+(z,x_0))=m, \lb{2.20} \\
\lim_{\varepsilon\downarrow 0} M_+(\lambda+
i\varepsilon,x_0) \text{ exists for a.e.\
$\lambda\in\bbR$}. \lb{2.21}
\end{gather}
Isolated poles of $M_+(z,x_0)$ and $-M_+(z,x_0)^{-1}$
are at most of first order, 
are real, and have a nonpositive residue.  \\
(ii)  $M_+(z,x_0)$ admits the representations
\begin{align}
M_+(z,x_0)&=F_+(x_0)+\int_\bbR
d\Omega_+(\lambda,x_0) \,
\big((\lambda-z)^{-1}-\lambda(1+\lambda^2)^{-1}\big)
\lb{2.22}  \\
&=\exp\bigg(C_+(x_0)+\int_\bbR d\lambda \, \Xi_+
(\lambda,x_0)
\big((\lambda-z)^{-1}-\lambda(1+\lambda^2)^{-1}\big)
\bigg), \lb{2.23}
\end{align}
where
\begin{align}
F_+(x_0)&=F_+(x_0)^*, \quad \int_\bbR \f{\Vert
d\Omega_+(\lambda,x_0)\Vert}{1+\lambda^2}<\infty,
\lb{2.24} \\
C_+(x_0)&=C_+(x_0)^*, \quad 0\le\Xi_+(\dott,x_0)
\le I_m \, \rm{  a.e.}
\lb{2.25}
\end{align}
Moreover,
\begin{align}
\Omega_+((\lambda,\mu],x_0)&=\lim_{\delta\downarrow
0}\lim_{\varepsilon\downarrow 0}\f1\pi
\int_{\lambda+\delta}^{\mu+\delta} d\nu \, \Im(
M_+(\nu+i\varepsilon,x_0)), \lb{2.26} \\
\Xi_+(\lambda,x_0)&=\lim_{\varepsilon\downarrow 0}
\f1\pi\Im(\ln(M_+(\lambda+i\varepsilon,x_0)) \text{ for
a.e.\ $\lambda\in\bbR$}.\lb{2.27}
\end{align}
(iii)  With $\phi(z,x,x_0)$ and $\theta(z,x,x_0)$ given in \eqref{2.4},
define the $m\times m$ matrices
\begin{equation}
u_+(z,x,x_0)=\theta(z,x,x_0) + \phi(z,x,x_0) M_+(z,x_0),
\quad x,x_0\in\bbR, \, z\in\bbC\backslash\bbR.
\lb{2.28}
\end{equation}
Then
\begin{equation}
u_+(z,\cdot,x_0)\in L^2((x_0,\infty))^{m\times m},
\lb{2.29}
\end{equation}
$u_+(z,x,x_0)$ is invertible, $u_+(z,x,x_0)$
satisfies the boundary condition of
$H_+$  at infinity (if any), and
\begin{equation}
\Im(M_+(z,x_0))=\Im(z)\int_{x_0}^{\infty}dx\,
u_+(z,x,x_0)^* u_+(z,x,x_0).\lb{2.30}
\end{equation}
\end{theorem}

\section{Atkinson's Leading-Order Argument}\lb{s3}

The purpose of this section is to prove the following result
due to Atkinson.

\begin{theorem} \mbox{\rm (Atkinson \cite{At88}.)} \lb{t3.1}
Assume Hypothesis~\ref{h2.1}
and denote by $C_\varepsilon\subset
\bbC_+$ the open sector with vertex at zero, symmetry
axis along the positive imaginary axis, and opening angle
$\varepsilon$, with $0<\varepsilon< \pi/2$. Let $M_+(z,x_0)$
be either the unique limit point or a point 
of the limit disk $D_+(z,x_0)$
associated with \eqref{2.2}. Then
\begin{equation}\lb{3.1}
M_+(z,x_0)\underset{\substack{\abs{z}\to\infty\\ z\in
C_\varepsilon}}{=} iz^{1/2}I_m + o(|z|^{1/2}).
\end{equation}
\end{theorem}

More precisely, Atkinson proved the analog of \eqref{3.1} in
\cite{At88}  for
more general Sturm-Liouville-type situations rather
than matrix-valued Schr\"odinger operators. To the best
of our knowledge, Atkinson's manuscript
\cite{At88}  remains unpublished. Since the result is
crucial for our principal Section~\ref{s4}, we will
specialize his arguments to the simplified situation at hand
and provide a complete account of his approach (at times going
beyond some of the arguments sketched in \cite{At88}).

One of our tools in studying $M_+(z,x_0)$ will be to relate it
to a matrix Riccati-type equation with solution $M_+(z,x)$
defined as follows. We first note that
\begin{equation}
M_+(z,x_0)=\ti u_+^\prime(z,x_0) \ti u_+(z,x_0)^{-1},
\quad z\in\bbC\backslash\bbR, \lb{3.2}
\end{equation}
where $\ti u_+(z,\dott)\in L^2([x_0,\infty))^{m\times m}$
is a nonnormalized solution of \eqref{2.2} satisfying
\begin{equation}
\ti u_+(z,x)=u_+(z,x,x_0)C_+,  \lb{3.3}
\end{equation}
with $u_+ $ defined in \eqref{2.28}, and with $C_+\in\bbC^{m\times m}$
a nonsingular matrix.  Varying the reference point $x_0$ we may define
\begin{equation}
M_+(z,x)=\ti u_+^\prime(z,x) \ti u_+(z,x)^{-1},
\quad z\in\bbC\backslash\bbR, \, x\geq x_0. \lb{3.4}
\end{equation}
Then $M_+(z,x)$ is independent of the chosen normalization
$C_+$ and is well-known (see, e.g.,
\cite{CGHL99}, \cite{GH97}, \cite{Jo87}) to satisfy
\eqref{3.5} below.

\begin{lemma}\lb{l3.2}
Assume Hypothesis~\ref{h2.1}, suppose
$z\in\bbC\backslash\bbR$,
$x\in [x_0,\infty)$, and define $M_+(z,x)$ as in
\eqref{3.4}. Then $M_+(z,x)$
satisfies  the standard Riccati-type equation,
\begin{equation}
M_+^\prime(z,x)+M_+(z,x)^2=Q(x)-z I_m. \lb{3.5}
\end{equation}
\end{lemma}
\begin{proof}
Differentiating \eqref{3.4} with respect to $x$ and
inserting
$\ti u_+^{''}(z,x)+(Q(x)-z)\ti u_+(z,x)=0$
immediately yields \eqref{3.5}.
\end{proof}

Utilizing
\begin{equation}
\Psi(z,x,x_1)=\Psi(z,x,x_0)\Psi(z,x_0,x_1) \text{ and }
\Psi(z,x_0,x_1)\Psi(z,x_1,x_0)=I_{2m},
\end{equation}
(cf.~\eqref{2.4}) one computes for all $x_1\in\bbR$
\begin{align}
&\theta(z,x,x_1)+\phi(z,x,x_1)M_+(z,x_1) \no \\
&=u_+(z,x,x_0)
[\theta(z,x_0,x_1)+\phi(z,x_0,x_1)M_+(z,x_1)], \quad
z\in\bbC\backslash\bbR. \lb{3.7}
\end{align}
Since $u_+(z,\cdot,x_0)\in L^2([R,\infty))^{m\times m}$
for all $R\in\bbR$, the left-hand side of \eqref{3.7} is
in $L^2([R,\infty))^{m\times m}$. Thus, for $M_+(z,x_1)$
defined according to \eqref{3.4}, one concludes that
\begin{equation}
M_+(z,x_1)\in D_+(z,x_1) \text{ for all }x_1\in\bbR
\end{equation}
since $M_+(z,x_0)\in D_+(z,x_0)$ by hypothesis. This
justifies our choice of notation $M_+(z,x)$ when compared
with the notation employed in Definition~\ref{d2.7}.

Atkinson's proof of Theorem~\ref{t3.1}, which we follow very
closely, is first concerned with the relationship between
$D(z,c,x_0)$,
described in Definition~\ref{d2.4},
and $D^\calR (z,c,x_0)$, a {\em Riccati disk}, which we
now define.

\begin{definition}\lb{d3.3}
With $Q(x)$ given in Hypothesis~\ref{h2.1} and
$z\in \bbC_+$,
$D^\calR (z,c,x_0)$  is defined to be the set of all
$M(z,x_0)\in\bbC^{m\times m}$ where
$\Im{(M(z,x_0))} > 0$, and such that the solution of
\begin{subequations}\lb{3.9}
\begin{align}\lb{3.9a}
\vT '(z,x)= \frac{1}{2}\begin{bmatrix} I_m+\vT (z,x)&
I_m-\vT (z,x)\end{bmatrix}
&\begin{bmatrix}-i|z|^{-1/2}(zI_m -Q(x))& 0\\0& i|z|^{1/2}I_m
\end{bmatrix}\times \no \\
&\times\begin{bmatrix} I_m+\vT (z,x)\\ I_m-\vT (z,x)\end{bmatrix},
\end{align}
\begin{align} \lb{3.9b}
\vT (z,x_0)&= (I_m+i|z|^{-1/2} M(z,x_0))(I_m-i|z|^{-1/2} M(z,x_0))^{-1}
\no \\
&= (i|z|^{1/2}I_m -M(z,x_0))(i|z|^{1/2}I_m +M(z,x_0))^{-1},
\end{align}
\end{subequations}
 satisfies
\begin{equation}\lb{3.10}
\vT (z,x)^* \vT (z,x)\le I_m  \text{ for all }  x\in[x_0,c].
\end{equation}
\end{definition}

Note that there is a correspondence between elements of 
$D^\calR (z,c,x_0)$
and certain solutions of \eqref{3.9a}. If $\vT (z,x) $ satisfies
\eqref{3.9a} and
\eqref{3.10} then note that a matrix $M(z,x_0)\in
D^\calR (z,c,x_0)$ is defined by the initial data $\vT (z,x_0) $
if $\vT (z,x_0)^*\vT (z,x_0)< I_m$.
  Note also that
\eqref{3.9a} can be written as
\begin{align}
\vT'=&\frac{1}{2}\begin{bmatrix}I_m+\vT & I_m-\vT\end{bmatrix}
\begin{bmatrix}I_m&0\\0&i|z|^{1/2}I_m\end{bmatrix}
(zA+B)
\begin{bmatrix}-i|z|^{-1/2} I_m& 0\\0& I_m \end{bmatrix}
\times \no \\
&\times \begin{bmatrix}I_m+\vT \\ I_m-\vT \end{bmatrix}.
\end{align}

\begin{remark} \lb{r3.4}
We emphasize that \eqref{3.9a}
formally results
from the Cayley-type transformation,
\begin{equation}
\vT (z,x)= (i|z|^{1/2}I_m -M(z,x))(i|z|^{1/2}I_m
+M(z,x))^{-1}, \lb{3.12}
\end{equation}
where $M(z,x)$ satisfies the Riccati-type equation \eqref{3.5}.
In the scalar context this corresponds to a conformal
mapping of the complex upper half-plane to the unit disk.  Moreover, we
note that Atkinson proves the asymptotic result of
Theorem~\ref{t3.1}
for elements of $D^\calR (z,c,x_0)$.  Recall that $D_+(z,x_0)\subseteq
D(z,c,x_0)$ (cf.~Remark~\ref{r2.6}). Theorem~\ref{t3.1} then
follows upon showing that $D(z,c,x_0)\subseteq D^\calR (z,c,x_0)$ \ \
for
$c>x_0$. This containment is shown in Theorem~\ref{t3.6}.
\end{remark}

\begin{lemma}\lb{l3.5}
If $u(z,x)$ is defined by \eqref{2.13} with $M\in
D(z,c,x_0)$, then $u(z,x) - i|z|^{-1/2} u'(z,x)$ is invertible for
all $z\in\bbC_+$, and all $x\in [x_0 , c]$.
\end{lemma}
\begin{proof}
With $M\in D(z,c,x_0)$ and with  $U(z,c)$ defined in
\eqref{2.13},  $U(z,c)$ satisfies \eqref{2.15}.
If $\Im{(\beta_2\beta_1^*)}> 0$, then $u'(z,c)= -\beta_2^{-1}\beta_1
u(z,c)$.  $U(z,c)$ has rank $m$, hence $u(z,c)$ is nonsingular. As a
result,
\begin{equation}
\Im{(u(z,c)^* u'(z,c))}= u(z,c)^*\beta_2^{-1}\Im{(\beta_2\beta_1^*)}
{\beta_2^{*}}^{-1} u(z,c)>0.
\end{equation}
However, if $\Im{(\beta_2\beta_1^*)}=0 $, then
\begin{equation}\lb{3.14}
0=\beta_1\beta_2^* - \beta_2\beta_1^* =
\begin{bmatrix}\beta_1& \beta_2
\end{bmatrix}\begin{bmatrix}\beta_2^*\\ -\beta_1^* \end{bmatrix}.
\end{equation}
By \eqref{2.5a}, $\rank [ \beta_1 \ \beta_2 ] =
\rank [ \beta_2 \ -\beta_1 ]^* = m$;  thus,  by \eqref{2.15}
and \eqref{3.14} there is a $w\in\bbC^{m\times m}$ such
that $U(z,c)= [
\beta_2 \ -\beta_1]^* w$.  From this, we obtain
\begin{equation}\lb{3.15}
\Im{(u(z,c)^*u'(z,c))}=w^*\Im{(\beta_2 \beta_1^*)}w=0.
\end{equation}
Hence for $\beta$ which satisfy \eqref{2.5},
\begin{equation}\lb{3.16}
\Im{(u(z,c)^*u'(z,c))\ge 0} .
\end{equation}
For all $x\in [x_0, c]$ and $z\in \bbC_+$, note that
\begin{equation}\lb{3.17}
(\Im{(u(z,x)^*u'(z,x))})' = - \Im{(z)} u(z,x)^*u(z,x) \le 0.
\end{equation}
Together, \eqref{3.16} and \eqref{3.17} imply that
\begin{equation}\lb{3.18}
\Im{(u(z,x)^*u'(z,x))\ge 0} \text{ for all }
x\in [x_0, c], \, z\in \bbC_+ .
\end{equation}

If for some $\xi \in [x_0, c] $ there  is an $f\in\bbC^m$,
$f\ne 0$, such that $u(z,\xi)f = i|z|^{-1/2} u'(z,\xi)f$,
then
\begin{equation}\lb{3.19}
-i|z|^{1/2}f^*u(z,\xi)^*u(z,\xi)f = f^*u(z,\xi)^*u'(z,\xi)f.
\end{equation}
Together, \eqref{3.18} and \eqref{3.19} imply that
$f^*u(z,\xi)^*u(z,\xi)f\le 0$, and hence that
$u(z,\xi)f=u'(z,\xi)f=0$.  Given the uniqueness of
solutions for
system \eqref{2.3}, we conclude that $0=u(z,x_0)f=f$;
thereby
producing a contradiction which completes the proof.
\end{proof}

\begin{theorem}\lb{t3.6}
$D(z,c,x_0)\subseteq D^\calR (z,c,x_0)$ \ for all $c>x_0$
and $z\in\bbC_+$.
\end{theorem}
\begin{proof}
With $z\in \bbC_+$ and $M=M(z,c,x_0,\beta) \in D(z,c,x_0)$,
we see that $\Im{(M)}>0 $ by \eqref{2.12}.  With
$u(z,x)$ defined
in \eqref{2.13}, then by Lemma~\ref{l3.5}, we can
define
\begin{equation}\lb{3.20}
\vT(z,x)= (u(z,x)+i|z|^{-1/2} u'(z,x))(u(z,x)-i|z|^{-1/2}
u'(z,x))^{-1}
\end{equation}
for $x\in [x_0,\ c]$.  As defined, $\vT (z,x)$ satisfies \eqref{3.9b}.
Then for $x\in [x_0\ c]$,  \eqref{3.18} implies that
\begin{equation}\lb{3.21}
\begin{split}
&2i|z|^{-1/2} \{u(z,x)^*u'(z,x) -{u^*}'(z,x)u(z,x) \} \\
&= -4|z|^{-1/2} \Im{(u^*(z,x)u'(z,x))}\le 0.
\end{split}
\end{equation}
This is equivalent to
\begin{equation}\lb{3.22}
\begin{split}
&(u(z,x)^* -i|z|^{-1/2} u'(z,x)^*) (u(z,x)
+i|z|^{-1/2} {u'}(z,x))  \\
&\le (u(z,x)^* +i|z|^{-1/2} u'(z,x)^*)(u(z,x)
-i|z|^{-1/2} {u'}(z,x)).
\end{split}
\end{equation}
Given the invertibility of $u(z,x) -i|z|^{-1/2}
{u'}(z,x)$ shown in Lemma~\ref{l3.5}, we infer
 that
$\vT(z,x)$ satisfies \eqref{3.10}. The proof of
this theorem will be completed upon showing that $\vT(z,x)$
satisfies \eqref{3.9a}.

First, observe that
\begin{equation}\lb{3.23}
\vT(z,x)(u(z,x)-i|z|^{-1/2} u'(z,x)) = u(z,x)
+ i|z|^{-1/2} u'(z,x),
\end{equation}
and hence that
\begin{equation}\lb{3.24}
I_m + \vT = 2u(u-i|z|^{-1/2} u')^{-1},\,\,\,
I_m - \vT = -2i|z|^{-1/2} u'(u-i|z|^{-1/2} u')^{-1}.
\end{equation}
Upon differentiating \eqref{3.23} and using the fact
that $u(z,x)$ satisfies \eqref{2.2}, we obtain
\begin{equation}
\vT'(u-i|z|^{-1/2} u') = (I_m- \vT)u' -i|z|^{-1/2}
(I_m+\vT)(zI_m -Q)u.
\end{equation}
By \eqref{3.24}, it follows that
\begin{equation}
\vT' = -\frac{1}{2i|z|^{-1/2}}(I_m -\vT)^2
- \frac{i|z|^{-1/2}}{2}(I_m +\vT)(zI_m -Q) (I_m +\vT),
\end{equation}
which is precisely \eqref{3.9a}.
\end{proof}
The following relation holds for the sets in
 Theorem~\ref{t3.6}, that is,
\begin{equation}
\overline{D(z,c,x_0)}=D^\calR (z,c,x_0) \text{ for all }c>x_0
\text{ and }z\in\bbC_+.
\end{equation}
This will be discussed in detail in \cite{CG00}.

An associated system for \eqref{3.9} is obtained by
introducing a change
of independent variable: With $z\in \bbC_\varepsilon$,\
$c\in [x_0 , \infty)$, and\ $x\in [x_0 , c]$, let
\begin{equation}\lb{3.28}
t\ = \ (x-x_0)|z|^{1/2}.
\end{equation}
Introducing
\begin{equation}
\vP (z,t) = \vT (z,x), \lb{3.29}
\end{equation}
with $x$ and $t$ related as in \eqref{3.28}, \eqref{3.9} becomes
\begin{subequations} \lb{3.30}
\begin{align}\lb{3.30a}
\vP '(z,t)= &\frac{1}{2}\begin{bmatrix} I_m+\vP (z,t)&
I_m-\vP (z,t)\end{bmatrix}
\begin{bmatrix}-i|z|^{-1}(zI_m -\hatt Q(t))& 0\\0& iI_m
\end{bmatrix} \times \no \\
&\times \begin{bmatrix} I_m+\vP (z,t)\\ I_m
-\vP (z,t)\end{bmatrix},
\end{align}
\begin{equation}\lb{3.30b}
\vP (z,0)=  (i|z|^{1/2}I_m -M(z,x_0))(i|z|^{1/2}I_m +M(z,x_0))^{-1},
\end{equation}
\end{subequations}
and \eqref{3.10} becomes
\begin{equation}\lb{3.31}
\vP (z,t)^*\vP (z,t)\le I_m  \text{ for all } t \in
[0,(c-x_0)|z|^{1/2}].
\end{equation}
Note that in \eqref{3.30a}
\begin{equation}\lb{3.32}
\hatt Q(t)=Q(x_0 + t|z|^{-1/2}).
\end{equation}
In \eqref{3.30} and \eqref{3.31}, we have a set
of conditions equivalent to \eqref{3.9} and \eqref{3.10}
for $D^\calR (z,c,x_0)$ given in Defintion~\ref{d3.3}.

Now consider a sequence, $z_n \in \bbC_{\varepsilon}$, such
that $|z_n| \to \infty$ as $n\to \infty$ and such that
\begin{equation}\lb{3.33}
0< \varepsilon <  \delta_n = \arg{(z_n)} < \pi - \varepsilon.
\end{equation}
By choosing an appropriate subsequence, we may assume that
\begin{equation}\lb{3.34}
\delta_n \to \delta \in [\varepsilon, \pi - \varepsilon].
\end{equation}
Let $\vP (z_n ,t)$ denote a corresponding sequence of functions
that satisfy \eqref{3.30a} and \eqref{3.31},
with initial
data, $\vP (z_n ,0)$, defined by \eqref{3.30b} for a
sequence of
points  $M(z_n,x_0)$, $n\in\bbN $, where each $M(z_n,x_0)$ is  chosen to
lie in the
disk $D^\calR (z_n,c,x_0)$. Note that as $z_n\to \infty$, the
intervals
described in
\eqref{3.31} eventually cover all compact
subintervals of
$ \bbR_+$. Given the uniform boundedness of
$\vP_n(t)=\vP (z_n ,t)$ described in \eqref{3.31},
we assume, upon passing to an appropriate subsequence
still denoted by  $\vP_n (0)$, that
\begin{equation}\lb{3.35}
\vP_n(0) = \vP (z_n,0) \rightarrow \vP_0(\delta), \text{ as }
n\rightarrow \infty,
\end{equation}
and as a consequence, that
\begin{equation}
{\vP_0(\delta)}^* \vP_0(\delta) \le I_m.
\end{equation}

With $\vP_0(\delta)$ defined in \eqref{3.35} as
$|z_n|\to\infty$, we
consider  a limiting system associated  with \eqref{3.30}:
\begin{subequations}\lb{3.37}
\begin{align}
\vE '(\delta,t)= &\frac{1}{2}
\begin{bmatrix}
I_m+ \vE (\delta,t) & I_m-\vE (\delta,t)
\end{bmatrix}
\begin{bmatrix}
-ie^{i\delta}I_m & 0\\ 0 & iI_m
\end{bmatrix}\times\no\\
&\times\begin{bmatrix}
I_m+\vE (\delta,t)\\ I_m- \vE (\delta,t)
\end{bmatrix},
\quad t\geq 0, \lb{3.37a}
\end{align}
\begin{equation}
\vE (\delta,0)= \vP_0(\delta).
\end{equation}
\end{subequations}

\begin{theorem}\lb{t3.7}
The solution $\vE(\delta,t) $ of \eqref{3.37}  satisfies
\begin{equation}\lb{3.38}
\vE (\delta,t)^* \vE(\delta,t) \le I_m
\end{equation}
for $ 0\le t < \infty$. Moreover, the
solutions $\vP _n (t)=\vP (z_n,t) $
of \eqref{3.30} converge to $\vE (\delta,t) $  uniformly on
$[0,T]$ for every $T>0$, as $n\to \infty $.
\end{theorem}
\begin{proof}
Let $T\in \bbR_+$ be the greatest value such that
\eqref{3.38} holds for
 $t\in [0,\ T] $.  We show that \eqref{3.38} must hold
for some $[0,\ T'] $ with  $ T' > T $, thus proving $T=\infty$.

The solution of  \eqref{3.37}, $\vE (\delta,t)$,
presumed
to be defined
 on  $[0,\ T]$, can be continued onto some $[0,\ T']$
with $T' > T$;
 $\vE (\delta,t)$ then satisfies
\begin{equation}\lb{3.39}
\vE (\delta,t)^* \vE (\delta,t) \le k^2 I_m
\end{equation}
for $0\le t \le T'$ and for some $k\ge 1$.

For brevity, let $\vP'_n = G_n(\vP_n,t) $ denote
\eqref{3.30a} with
$z=z_n$, and let $\vE'= H(\vE,t) $ denote \eqref{3.37a}
in
the following. Integrating \eqref{3.37} and
\eqref{3.30}, we obtain
\begin{align}
\vP_n (t) -\vE (\delta,t) &= \vP_n(0) -\vP_0(\delta)
+ \int_0^t  ds \{
G_n(\vE,s) - H(\vE,s)\}  \no \\
&\quad + \int_0^t ds \{ G_n(\vP_n,s) - G_n(\vE,s)\}. \lb{3.40}
\end{align}
Note that
\begin{align}
G_n(\vE ,s) - H(\vE ,s) &= \frac{1}{2}i(e^{i\delta} -
e^{i\delta_n}) (I_m +\vE (\delta,s))^2  \no \\
&\quad + \frac{1}{2}i|z_n|^{-1}(I_m
+ \vE (\delta,s))\hatt Q(s)(I_m + \vE (\delta,s)),
\end{align}
and by  \eqref{3.39}, as $n\to\infty$ and for $t \in [0,T']$, that
\begin{equation}
\frac{1}{2}i(e^{i\delta} - e^{i\delta_n}) \int_0^t ds\ (I_m + \vE
(\delta,s) )^2 = O (\delta - \delta_n),
\end{equation}
and by \eqref{3.28}, \eqref{3.32}, and \eqref{3.39}, that
\begin{align}
&\frac{1}{2}i|z_n|^{-1}\int_0^t ds\ (I_m + \vE (\delta,s))
\hatt Q(s)(I_m + \vE (\delta,s))  \no \\
&= |z_n|^{-1/2}O \left ( \int_{x_0}^{x_0
+t|z_n|^{-1/2}}dx\ || Q(x)
||  \right )   \lb{3.43} \\
&= o (|z_n|^{-1/2}) \text{ as } n\to\infty. \lb{3.44}
\end{align}
Together with \eqref{3.35}, we observe that as $n \to \infty$
\eqref{3.40} is equivalent to
\begin{equation}
\vP_n (t) -\vE (\delta,t) = o (1) +  \int_0^t ds
(G_n(\vP_n,s) - G_n(\vE,s))
\end{equation}
uniformly so, for $0\le t \le T' $.
Since $|| \vP _n || \le 1 $ and $|| \vE || \le k$
for $ t\in [0, \ T' ] $,
\begin{equation}
|| G_n ( \vP _n,s) - G_n( \vE , s) || \le
\frac{1}{2}(3 + k)|| \vP _n (s)-\vE (\delta,s)||
(2 + |z_n|^{-1}\ ||\hatt Q (s)||).
\end{equation}
Thus \eqref{3.40} yields
\begin{equation}\lb{3.47}
\begin{split}
&|| \vP_n (t) -\vE (\delta,t)  ||  \\
&\le o(1) + \int_0^t ds\
 \frac{1}{2}(3 + k)|| \vP _n (s)-\vE (\delta,s) ||
 (2 + |z_n|^{-1}\ ||\hatt Q(s) ||).
\end{split}
\end{equation}
In light of \eqref{3.44}, an application of Gronwall's
inequality to \eqref{3.47} yields
\begin{equation}
\vP_n(t)  - \vE (\delta,t)  \to 0 \text{ as } n\to\infty
\lb{3.48}
\end{equation}
uniformly for $0\le t\le T'$.  For $n$ sufficiently
large and $t\in [0, \ T']$, $\vP_n (t) $ satisfies
\eqref{3.31}
with $z=z_n$, and $ \vE (\delta,t) $ satisfies \eqref{3.39},
hence $ \vE (\delta,t) $ satisfies \eqref{3.38} for
$t\in [0,\ T']$, where $T'>T$.
\end{proof}

The following computation identifies $\vP_0 (\delta)$
and proves
that $\vE (\delta,t)$ is constant with respect to $t\geq 0$.

\begin{corollary} \lb{c3.8}
\begin{equation} \lb{3.49}
\vE
(\delta,t)=\vE(\delta,0)=\vP_0(\delta)
=\frac{1-\exp(i\delta/2)}{1+\exp(i\delta/2)}I_m,
\quad t\geq 0.
\end{equation}
\end{corollary}
\begin{proof}
Utilizing the standard connection between the explicit
exponential solutions of the  second-order Schr\"odinger equation
\eqref{2.2a} and the Riccati equation \eqref{3.5} (in
the special case
$Q(x)=0$, cf.~\eqref{3.2} and Lemma~\ref{l3.2}), performing a
conformal map of the type
\eqref{3.12}, and the  variable transformations \eqref{3.28} and
\eqref{3.29} then yields the following solution for
\eqref{3.30a},
\begin{align}
\vP (z,t)=&\big(-(i|z|^{1/2}-iz^{1/2})(M_0(z)+iz^{1/2}I_m)+
\exp(-2it\exp(i\delta/2))\times \no \\
&\times (i|z|^{1/2}+iz^{1/2})(M_0(z)-iz^{1/2}I_m)\big)\times
\no \\
&\times \big(-(i|z|^{1/2}+iz^{1/2})(M_0(z)+iz^{1/2}I_m)+
\exp(-2it\exp(i\delta/2))\times \no \\
&\times (i|z|^{1/2}-iz^{1/2})
(M_0(z)-iz^{1/2}I_m)\big)^{-1}, \lb{3.50}
\end{align}
associated with the general initial condition
\begin{equation}
\vP (z,0)=\big(i|z|^{1/2}I_m -M_0(z)\big)
\big(i|z|^{1/2}I_m +M_0(z) \big)^{-1} \lb{3.51}
\end{equation}
for some $M_0(z)\in\bbC^{m\times m}$ with $\Im (M_0(z))>0$,
$z\in\bbC_+$.
Since by hypothesis $0<\delta<\pi$, the exponential terms in
\eqref{3.50} enforce
\begin{equation}
\|\vP (z,t)\| > 1 \text{ as } t\uparrow\infty
\end{equation}
unless
\begin{equation}
M_0(z)=iz^{1/2}I_m,
\end{equation}
implying \eqref{3.49}.
\end{proof}

Given these facts we  proceed to the

\vspace*{1mm}

\noindent {\it Conclusion of the proof of Theorem~\ref{t3.1}.}

\vspace*{1mm}

With $M(z_n,x_0)\in D^{\calR}(z,c,x_0)$,  and for $\delta$ as in
\eqref{3.34} let
\begin{equation}
C(\delta)=\frac{1-\exp(i\delta/2)}{1+\exp(i\delta/2)}. \lb{3.54}
\end{equation}
 From \eqref{3.30b}, \eqref{3.35},\eqref{3.48},  and
\eqref{3.49}, we infer that
\begin{equation}
(I_m +i|z_n|^{-1/2} M(z_n,x_0))(I_m-i|z_n|^{-1/2} M(z_n,x_0))^{-1}
\to C(\delta) I_m \text{ as } n\to \infty.
\end{equation}
As a consequence,
\begin{align}
i|z_n|^{-1/2} M(z_n,x_0) &= \frac{C(\delta)-1
+o(1)}{C(\delta)+1 +o(1)}I_m=
\frac{C(\delta)-1}{C(\delta)+1}(1+o(1))I_m \no \\
&=-e^{i\delta/2}(1+o(1))I_m. \lb{3.56}
\end{align}
Thus, by \eqref{3.54} and \eqref{3.56} we conclude that
\begin{equation}\lb{3.57}
M(z_n,x_0)= i z_n^{1/2}(1+o(1)) I_m,
\end{equation}
and hence by Theorem~\ref{t3.6}, that \eqref{3.1} holds.
\hspace*{5.33cm} $\square$\\

In \eqref{3.1} an asymptotic expansion is given that is uniform with
respect to $\arg(z)$ for $|z| \to \infty$ in $C_{\varepsilon}$.  However,
we observe that the proof just completed shows more:  Allowing the
reference point $x_0$ to vary, the asymptotic expansion given in
\eqref{3.1} is also uniform in $x_0$ for $x_0$ in a compact subset of
$\bbR$.

\begin{theorem} \lb{t3.9}
Assume Hypothesis~\ref{h2.1}, let $z\in\bbC_+$,
$x\in \bbR$, and denote by $C_\varepsilon\subset
\bbC_+$ the open sector with vertex at zero, symmetry
axis along the positive imaginary axis, and opening angle
$\varepsilon$, with $0<\varepsilon< \pi/2$. Let $M_+(z,x)$,
$x\geq x_0$,
be as in Definition~\ref{d2.7}. Then
\begin{equation}
M_+(z,x)\underset{\substack{\abs{z}\to\infty\\ z\in
C_\varepsilon}}{=} iz^{1/2}I_m + o(|z|^{1/2}) \lb{3.58}
\end{equation}
uniformly with respect to $\arg\,(z)$ for $|z|
\to \infty$ in $C_\varepsilon$ and uniformly in $x$ as long as $x$
varies in compact subsets of $[x_0,\infty)$.
\end{theorem}
\begin{proof}
Though stated for elements of the Weyl disk $D_+(z,x_0)$, the 
proof of
Theorem~\ref{t3.1} shows that the asymptotic expansion given 
in \eqref{3.1}
holds, uniformly with respect to $\arg\,(z)$ for $|z|
\to \infty$ in $C_\varepsilon$, for all elements of the Riccati disk
$D^{\calR}(z,c,x_0)$ as noted in \eqref{3.57}.

Note that the system \eqref{3.37} is independent of the reference
point $x_0$. Recall that $\delta$, defined in \eqref{3.34},  is
determined by our apriori choice of the sequence $z_n$ subject only to
$z_n$  being in $C_\varepsilon$ (c.f. \eqref{3.33}). Next note that
$\vP_0(\delta)$, which is defined as a limit in \eqref{3.35}, which is
described explicitly in Corollary~\ref{c3.8}, and which gives 
solutions of
\eqref{3.37} which satisfy \eqref{3.38} for $0\le t<\infty$, is
also independent of the reference point $x_0$. Thus, had we chosen
$x_0'\ne x_0$ as our reference point at the start, the asymptotic
analysis begun in Theorem~\ref{t3.7} and continued in
\eqref{3.54}--\eqref{3.57} would remain the same after the variable
change in \eqref{3.28} except for the integral expression present in
\eqref{3.43} in which $x_0$ would be replaced by $x_0'$.  However,
given the local integrability assumption on
$Q(x)$ present in Hypothesis~\ref{h2.1}, we see that the integral
expression in \eqref{3.43} is uniformly continuous for $x_0$ in 
a compact
subset of $\bbR$.  Thus \eqref{3.44} and consequently \eqref{3.48} are
uniform for $t$, and for $x_0$, in compact subsets of $\bbR$.

As a consequence, we see that \eqref{3.57} holds for elements of the
Riccati disk $D^{\calR}(z,c,x_0)$, that this asymptotic expansion is
uniform with respect to $\arg\,(z)$ for $|z|\to \infty$ in 
$C_\varepsilon$
and that it is uniform in $x_0$ as long as $x_0$ varies in 
compact subsets
of $\bbR$. The asymptotic expansion described in \eqref{3.58} 
then follows by Theorem~\ref{t3.6}
\end{proof}

We will see in the next section that the remainder term
$o(|z|^{1/2})$ in \eqref{3.58} can be improved to $o(1)$.

\section{Higher-Order Asymptotic Expansions}\lb{s4}

In this section we shall prove our principal result, the
asymptotic high-energy expansion of $M_+(z,x)$ to
arbitrarily high orders in sectors of the type
$C_\varepsilon\subset\bbC_+$ as defined in
Theorem~\ref{t3.1}.

Throughout this section we choose $\Im(z^{1/2})>0$ for
$z\in\bbC_+$. We also recall the following notion: $x\in [a,b)$ (resp.,
$x\in (a,b]$) is called a  right (resp., left) Lebesgue point of an
element
$q\in L^1 ((a,b))$, $a<b$, if $\int_0^\varepsilon
dx^\prime \,  |q(x+x^\prime)-q(x)|=o(\varepsilon)$ (resp., 
$\int_0^\varepsilon dx^\prime \, |q(x-x^\prime)-q(x)|=o(\varepsilon)$) 
as $\varepsilon\downarrow 0$. Similarly, $x\in (a,b)$ is called
a Lebesgue point of $q\in L^1 ((a,b))$ if
$\int_{-\varepsilon}^\varepsilon dx^\prime
\,  |q(x+x^\prime)-q(x)|=o(\varepsilon)$ as $\varepsilon\downarrow 0$. 
The set of all such points is then denoted by the right (resp., left) 
Lebesgue set of $q$ on $[a,b]$ in the former case and simply the Lebesgue
set of $q$ on $[a,b]$ in the latter case. The
analogous notions are applied to $m\times m$ matrices $Q\in L^1
((a,b))^{m\times m}$ by simultaneously considering all $m^2$ entries of
$Q$. The right (resp., left) Lebesgue set of $Q$ on $[a,b]$ is then
simply the intersection of the right (resp., left) Lebesgue sets of
$Q_{j,k}$ for all $1\leq j,k\leq m$, and similarly for the Lebesgue 
set of $Q$, etc.
 
We start with the simple case where $Q$ has compact support
and hence the underlying half-line matrix-valued Schr\"odinger
operator is in the limit point case at $+\infty$. 

\begin{lemma}\lb{l4.1}
Fix $x_0\in\bbR$ and let $x\geq x_0$. In addition to Hypothesis~\ref{h2.1}
suppose that $\widetilde Q$ has compact support in $[x_0,\infty)$, that
$\widetilde Q^{(N-1)}\in L^1([x_0,\infty))^{m\times m}$ for some
$N\in\bbN$, and that $x$ is a right Lebesgue point of $\widetilde
Q^{(N-1)}$. Denote by $\widetilde M_+(z,x)$  the Weyl-Titchmarsh matrix
associated with $\widetilde Q$. Then, as $\abs{z}\to\infty$ in
$C_\varepsilon$, $\widetilde M_+(z,x)$ has an asymptotic
expansion of the form $(\Im(z^{1/2})>0$, $z\in\bbC_+)$
\begin{equation}
\widetilde M_+(z,x)\underset{\substack{\abs{z}\to\infty\\ z\in
C_\varepsilon}}{=}
 i I_m z^{1/2}+\sum_{k=1}^N \tilde m_{+,k}(x)z^{-k/2}+
o(|z|^{-N/2}), \quad N\in\bbN. \lb{4.1}
\end{equation}
The expansion \eqref{4.1} is uniform with respect to
$\arg\,(z)$ for $|z|
\to \infty$ in
$C_\varepsilon$ and uniform in $x$ as long as $x$ varies in compact
subintervals of $[x_0,\infty)$ intersected with  the right Lebesgue set
of $\widetilde Q^{(N-1)}$.  The expansion  coefficients $\tilde
m_{+,k}(x)$ can be recursively computed from
\begin{align}
\tilde m_{+,1}(x)&=\f1{2i} \widetilde Q(x),
\quad \tilde m_{+,2}(x)= \f1{4} \widetilde Q^\prime(x),
\no \\
\tilde m_{+,k+1}(x)&=\f{i}2\bigg(\tilde m_{+,k}^\prime(x)+
\sum_{\ell=1}^{k-1}\tilde m_{+,\ell}(x)
\tilde m_{+,k-\ell}(x) \bigg),
\quad k\ge 2. \lb{4.2}
\end{align}

\noindent If one merely assumes 
$\wti Q\in L^1([x_0,\infty))^{m\times m}$ with compact support in 
$[x_0,\infty)$,  then
\eqref{4.1}  holds with $N=0$ {\rm (}interpreting $\sum_{k=1}^0 \cdot
=0${\rm )}, uniformly with respect to $\arg (z)$ for $|z|\to \infty$ in
$C_\varepsilon$ and uniformly in $x$ as long as $x$ varies in compact
subsets of $[x_0,\infty)$.
\end{lemma}
\begin{proof}
In the following let $z\in\bbC_+$, $\Im(z^{1/2}) > 0$, and $x\geq x_0$.
The existence of an expansion of the type \eqref{4.1} is shown as
follows.  First one considers a matrix Volterra integral equation of the
type (cf.~\cite[Ch. I]{AM63}, \cite{Co77}, \cite{NJ55},
\cite{Sc83}, \cite{WK74})
\begin{equation}
\tilde u_+(z,x)=\exp(iz^{1/2}x)I_m - \int_x^\infty dx'\,
\f{\sin(z^{1/2}(x-x'))}{z^{1/2}}\widetilde Q(x')\tilde u_+(z,x'),  
\lb{4.3} 
\end{equation}
and observes that the solution
$\tilde u_+(z,x)$ of \eqref{4.3} satisfies
$\tilde u_+(z,\cdot)\in L^2([x_0,\infty))^{m\times m}$ in
accordance with \eqref{3.3} and \eqref{3.4}. Next, introducing
\begin{equation}
\tilde v_+(z,x)=\exp(-iz^{1/2}x)\tilde u_+(z,x), \lb{4.3a}
\end{equation}
one rewrites \eqref{4.3} in the form
\begin{equation}
\tilde v_+(z,x)=I_m - \f{1}{2iz^{1/2}}\int_x^\infty dx'\,
[1-\exp(-2iz^{1/2}(x-x'))]\widetilde Q(x')\tilde v_+(z,x'),
\lb{4.3b} 
\end{equation}
and thus infers, 
\begin{equation}
\widetilde M_+(z,x)=\tilde u_+'(z,x)\tilde u_+(z,x)^{-1} 
=iz^{1/2}I_m + \tilde v_+'(z,x)\tilde v_+(z,x)^{-1}. \lb{4.4}
\end{equation}
Iterating \eqref{4.3b} then yields
\begin{equation}
\|\tilde v_+(z,x) \| \leq C, \quad z\in\bbC_+, \,\, \Im(z^{1/2}) > 0, 
\,\, x\geq x_0 \lb{4.4c}
\end{equation}
for some $C>0$ depending on $\widetilde Q$.
Finally, we need one more ingredient, proven in \cite[Lemma~3]{Ry99} 
using appropriate maximal functions.  Let $q\in L^1 ([x_0,\infty))$, 
$\supp(q)\subseteq [x_0,R]$, for some $R>0$, and suppose $x\in [x_0,R]$
is a right Lebesgue point of $q$. Then
\begin{equation}
\int_x^R dx^\prime \, \exp(2iz^{1/2}(x^\prime -x))q(x^\prime)
\underset{\substack{\abs{z}\to\infty\\ z\in
C_\varepsilon}}{=}-\f{q(x)}{2iz^{1/2}} + o(|z|^{-1/2}). \lb{4.3d} 
\end{equation}
An alternative proof of \eqref{4.3d} follows from 
\cite[Theorem~I.13]{Ti86}, which implies
\begin{equation}
\lim_{\substack{\abs{z}\to\infty\\ z\in
C_\varepsilon}}z^{-1/2} \int_x^R dx^\prime \, \exp(2iz^{1/2}(x^\prime
-x)) |q(x^\prime) -q(x)| = 0 \lb{4.3e}
\end{equation}
for any right Lebesgue point $x$ of $q$.
Given these facts, one iterates \eqref{4.3b} and its $x$-derivative,
integrates by parts, applies \eqref{4.3d} to $q=\|Q_{j,k}\|$ for all 
$1\leq j,k\leq m$, and estimates  $\|\tilde v_+(z,x^\prime)\|$ by
\eqref{4.4c}.  Inserting the expansions for $\tilde v_+'(z,x)$ and 
$\tilde v_+(z,x)^{-1}$into \eqref{4.4}  (using a geometric series
expansion for $\tilde v_+(z,x)^{-1}$) then yields the existence of an
expansion of the type \eqref{4.1}. The actual expansion coefficients and
the associated recursion relation \eqref{4.2} then follow upon inserting
expansion \eqref{4.1} into the Riccati-type equation \eqref{3.5}. The
assertion following \eqref{4.2} is an immediate consequence of
\eqref{4.3b} and its derivative with respect to $x$, \eqref{4.4}, 
and the Riemann-Lebesgue lemma.
\end{proof}

The corresponding scalar case $m=1$ (for a sufficiently smooth
coefficient $q$) was treated in this manner in \cite{GHSZ95}. The 
original version of our preprint assumed 
$\widetilde Q^{(N)}\in L^1([x_0,\infty))^{m\times m}$ in 
Lemma~\ref{l4.1}. Prompted by a recent 
preprint by A.~Rybkin \cite{Ry99}, who  
used an $m$-function-type approach to trace formulas for scalar 
Schr\"odinger operators and succeeded in removing any smoothness
hypotheses  on $q$, we reconsidered our original approach and extended
the  asymptotic expansion \eqref{4.1} to the case 
$\widetilde Q^{(N-1)}\in L^1 (\bbR)$ and $x$ a right Lebesgue point of 
$\widetilde Q^{(N-1)}$.

Next we recall an elementary result on finite-dimensional
evolution equations essentially taken from \cite{MPS90}.
\begin{lemma} \mbox{\rm (\cite{MPS90}.)} \lb{l4.2}
Let $A\in L^1_{\loc}(\bbR)^{m\times m}$. Then any
$m\times m$ matrix-valued solution $X(x)$ of
\begin{equation}
X'(x)=A(x)X(x)+X(x)A(x) \text{ for ~a.e. } x\in\bbR, \lb{4.5}
\end{equation}
is of the type
\begin{equation}
X(x)=Y(x)CZ(x), \lb{4.6}
\end{equation}
where $C$ is a constant $m\times m$ matrix and $Y(x)$ is a
fundamental system of solutions of
\begin{equation}
\Psi'(x)=A(x)\Psi(x) \lb{4.7}
\end{equation}
and $Z(x)$ is a fundamental system of solutions of
\begin{equation}
\Phi'(x)=\Phi(x)A(x). \lb{4.8}
\end{equation}
\end{lemma}
\begin{proof}
Clearly \eqref{4.6} satisfies \eqref{4.5} since
\begin{equation}
X'=Y'CZ+YCZ'=AYCZ+YCZA=AX+XA. \lb{4.9}
\end{equation}
Conversely, let $X$ be a solution of \eqref{4.5} and $Y$ a
fundamental matrix of solutions of \eqref{4.7}. Define
\begin{equation}
K(x)=Y(x)^{-1}X(x), \text{ that is, } X(x)=Y(x)K(x).
\lb{4.10}
\end{equation}
Then
\begin{equation}
X'=Y'K+YK'=AYK+YK'=AX+YK'=AX+XA \lb{4.11}
\end{equation}
implies
\begin{equation}
YK'=XA, \quad K'=Y^{-1}XA=KA. \lb{4.12}
\end{equation}
Thus, there exists a constant $m\times m$ matrix $C$
(possibly singular), such that
\begin{equation}
K(x)=CZ(x), \lb{4.13}
\end{equation}
with $Z$ a fundamental matrix of solutions of \eqref{4.8}.
Hence, $X=YCZ$.
\end{proof}

The next result provides the proper extension of
Proposition~2.1  in \cite{GS98} to the matrix-valued case.

\begin{lemma} \lb{l4.3}
Suppose $Q_1,Q_2 \in L^1_{\loc} (\bbR)^{m\times m}$ with
$Q_1(x)=Q_2(x)$ for~a.e. $x\in [x_0,x_1]$, $x_0<x_1$.
Denote by $M_{j,+}(z,x)$ the Weyl-Titchmarsh
matrices corresponding to $Q_j$, $j=1,2$. Then
\begin{align}
&(M_{1,+}'(z,x)-M_{2,+}'(z,x)) \no \\
&=-(1/2)(M_{1,+}(z,x)+M_{2,+}(z,x))
(M_{1,+}(z,x)-M_{2,+}(z,x)) \lb{4.14} \\
&\quad -(1/2)(M_{1,+}(z,x)-M_{2,+}(z,x))(M_{1,+}(z,x)
+M_{2,+}(z,x))  \no
\end{align}
for~a.e. $x\in [x_0,x_1]$.
\end{lemma}
\begin{proof}
This is obvious from \eqref{3.5}.
\end{proof}
\begin{lemma} \lb{l4.4}
Suppose $Q_1,Q_2 \in L^1_{\loc} (\bbR)^{m\times m}$ and
denote by $M_{j,+}(z,x)$ the corresponding Weyl-Titchmarsh
matrices associated with $Q_j$, $j=1,2$. Define
\begin{equation}
A_+(z,x)=-(1/2)(M_{1,+}(z,x)+M_{2,+}(z,x)) \lb{4.15}
\end{equation}
for $x\in [x_0,x_1]$, $x_0<x_1$ and assume $Y_+(z,x)$
and $Z_+(z,x)$ to be fundamental matrix solutions
of
\begin{equation}
\Psi'(z,x)=A_+(z,x)\Psi(z,x) \text{ and }
\Phi'(z,x)=\Phi(z,x)A_+(z,x) \lb{4.16}
\end{equation}
respectively, with
\begin{equation}
Y_+(z,x_1)=Z_+(z,x_1)=I_m. \lb{4.17}
\end{equation}
Then, as $|z|\to\infty$, $z\in C_\varepsilon$,
\begin{equation}
\|Y_+(z,x_0)\|, \|Z_+(z,x_0)\| \leq
\exp(-(x_1-x_0)\Im(z^{1/2})(1+o(1))). \lb{4.18}
\end{equation}
\end{lemma}
\begin{proof}
Define $B_+(z,x)$ by
\begin{equation}
B_+(z,x)=A_+(z,x)+iz^{1/2}I_m, \lb{4.19}
\end{equation}
then
\begin{equation}
\int_{x_0}^{x_1} dx\, \|B_+(z,x)\|
\underset{\substack{\abs{z}\to\infty\\ z\in
C_\varepsilon}}{=}  o(|z|^{1/2}) \lb{4.20}
\end{equation}
due to the uniform nature of the asymptotic expansion
\eqref{3.58} for $x$ varying in compact intervals. Next,
introduce
\begin{equation}
E_+(z,x,x_1)=\exp(i(x_1-x)z^{1/2})I_m, \quad x\leq x_1,
\lb{4.21}
\end{equation}
then
\begin{align}
Y_+(z,x)&=E_+(z,x,x_1)-\int_x^{x_1}dx'\,E_+(z,x,x')
B_+(z,x')Y_+(z,x'), \lb{4.22} \\
Z_+(z,x)&=E_+(z,x,x_1)-\int_x^{x_1}dx'\,Z_+(z,x')
B_+(z,x')E_+(z,x,x'). \lb{4.23}
\end{align}
Using
\begin{equation}
\|E_+(z,x_0,x_1)\|\leq\exp(-(x_1-x_0)\Im(z^{1/2})),
\lb{4.24}
\end{equation}
a standard Volterra-type iteration argument in
\eqref{4.22}, \eqref{4.23} then yields
\begin{equation}
\|Y_+(z,x_0)\|, \|Z_+(z,x_0)\| \leq
\exp\left(-(x_1-x_0)\Im(z^{1/2})+\int_{x_0}^{x_1}dx
\,\|B_+(z,x)\|\right), \lb{4.25}
\end{equation}
and hence \eqref{4.18}.
\end{proof}
\begin{theorem} \lb{t4.5}
Suppose $Q_1,Q_2 \in L^1_{\loc} (\bbR)^{m\times m}$ with
$Q_1(x)=Q_2(x)$ for~a.e. $x\in [x_0,x_1]$, $x_0<x_1$ and
let $M_{j,+}(z,x)$ be the corresponding Weyl-Titchmarsh
matrices associated with $Q_j$, $j=1,2$. Then, as
$\abs{z}\to\infty$ in $C_\varepsilon$,
\begin{equation}
\|M_{1,+}(z,x_0)-M_{2,+}(z,x_0)\|\leq
C(1+|z|^{1/2})\exp(-2(x_1-x_0)\Im(z^{1/2})(1+o(1)))
\lb{4.26}
\end{equation}
for some constant $C>0$.
\end{theorem}
\begin{proof}
Define for $x\in [x_0,x_1]$, $z\in\bbC\backslash\bbR$
\begin{align}
X_+(z,x)&=M_{1,+}(z,x)-M_{2,+}(z,x),  \\
A_+(z,x)&=-(1/2)(M_{1,+}(z,x)+M_{2,+}(z,x)).
\end{align}
By Lemma~\ref{l4.3},
\begin{equation}
X_+'=A_+X_++X_+A_+
\end{equation}
and hence by Lemma~\ref{l4.2},
\begin{equation}
X_+(z,x)=Y_+(z,x)X_+(z,x_1)Z_+(z,x), \lb{4.30}
\end{equation}
where $Y_+(z,x)$ and $Z_+(z,x)$ are fundamental matrices of
\begin{equation}
\Psi'(z,x)=A_+(z,x)\Psi(z,x) \text{ and }
\Phi'(z,x)=\Phi(z,x)A_+(z,x),
\end{equation}
respectively, with
\begin{equation}
Y_+(z,x_1)=I_m, \quad Z_+(z,x_1)=I_m.
\end{equation}
By Lemma~\ref{l4.4},
\begin{equation}
\|Y_+(z,x_0)\|, \|Z_+(z,x_0)\| \leq
\exp(-(x_1-x_0)\Im(z^{1/2})(1+o(1))) \lb{4.33}
\end{equation}
as $|z|\to\infty$, $z\in C_\varepsilon$. Thus, as
$|z|\to\infty$, $z\in C_\varepsilon$,
\begin{align}
\|X_+(z,x_0)\|&\leq
\|X_+(z,x_1)\|\,\|Y_+(z,x_0)\|\,\|Z_+(z,x_0)\| \no \\
&\leq C(1+|z|^{1/2})\exp(-2(x_1-x_0)\Im(z^{1/2})(1+o(1)))
\lb{4.34}
\end{align}
for some constant $C>0$ by \eqref{3.58}, \eqref{4.30}, and
\eqref{4.33}.
\end{proof}
The result \eqref{4.26} can be slightly improved as will be
discussed in Remark~\ref{r4.6a}.

Given these preparations we can now drop the compact
support assumption on $Q$ in Lemma~\ref{l4.1} and hence arrive at
the principal result of this paper.

\begin{theorem} \lb{t4.6}
Fix $x_0\in\bbR$. In addition to Hypothesis~\ref{h2.1}
suppose that for some $N\in\bbN$,
$Q^{(N-1)}\in L^1([x_0,c])^{m\times m}$ for all $c>x_0$, and that $x_0$ is
a right Lebesgue point of $Q^{(N-1)}$. Let $M_+(z,x_0)$
be either the unique limit point or a point of the limit disk
$D_+(z,x_0)$ associated with \eqref{2.2}. Then, as
$\abs{z}\to\infty$ in $C_\varepsilon$, $M_+(z,x_0)$ has an asymptotic
expansion of the form $(\Im(z^{1/2})>0$, $z\in\bbC_+)$
\begin{equation}
M_+(z,x_0)\underset{\substack{\abs{z}\to\infty\\ z\in
C_\varepsilon}}{=} i I_m z^{1/2}+\sum_{k=1}^N m_{+,k}(x_0)z^{-k/2}+
o(|z|^{-N/2}), \quad N\in\bbN. \lb{4.35}
\end{equation}
The expansion \eqref{4.35} is uniform with respect to $\arg\,(z)$ for $|z|
\to \infty$ in $C_\varepsilon$. The expansion  coefficients $m_{+,k}(x_0)$
can be recursively computed from \eqref{4.2} {\rm (}replacing $\tilde
m_{+,k}(x)$ by $m_{+,k}(x)${\rm )}.

\noindent If one merely assumes Hypothesis~\ref{h2.1}, then \eqref{4.35} 
holds with $N=0$ {\rm (}interpreting $\sum_{k=1}^0 \cdot =0${\rm )},
uniformly with respect to $\arg (z)$ for $|z|\to \infty$ in
$C_\varepsilon$.
\end{theorem}
\begin{proof}
Define
\begin{equation}
\widetilde Q(x)=\begin{cases} Q(x) &\text{ for }
x\in [x_0,x_1], \quad x_0<x_1 \no \\
0 \hspace*{6mm}&\text{ otherwise} \lb{4.36} \end{cases}
\end{equation}
and apply Theorem~\ref{t4.5} with $Q_1=Q$, $Q_2=\widetilde
Q$. Then (in obvious notation)
\begin{equation}
\|M_+(z,x_0)-\widetilde M_+(z,x_0)\|\leq C(1+|z|^{1/2})
\exp(-2(x_1-x_0)\Im(z^{1/2})(1+o(1)))
\end{equation}
as $|z|\to\infty$, $z\in C_\varepsilon$, and hence the
asymptotic expansion \eqref{4.1} for $\widetilde M_+(z,x_0)$
in Lemma~\ref{l4.1} coincides with that of $M_+(z,x_0)$ (this also  
applies to the case $N=0$).
\end{proof}

\begin{remark} \lb{r4.6a}
A comparison of \eqref{3.58} and \eqref{4.35} with $N=0$ shows that one
can improve the remainder term $o(|z|^{1/2})$ in \eqref{3.58}
to $o(1)$. In particular, this also yields an improvement of our rough
estimate \eqref{4.34} to the effect that the term $(1+|z|^{1/2})$
can be replaced by $1$ in \eqref{4.34} and hence in \eqref{4.26}.
\end{remark}

In analogy to Theorem~\ref{t3.9}, the asymptotic expansion
\eqref{4.35} extends to one for $M_+(z,x)$ valid uniformly
with respect to $x$ as long as $x$ varies in compact subintervals
of $[x_0,\infty)$ intersected with the right Lebesgue set of 
$Q^{(N-1)}$. Since this result represents a fundamental input for 
deriving trace formulas in \cite{CGHL99}, we provide the facts in 
the next theorem.

\begin{theorem} \lb{t4.7}
Fix $x_0\in\bbR$ and let $x\geq x_0$. In addition to
Hypothesis~\ref{h2.1} suppose that for some $N\in\bbN$,
$Q^{(N-1)}\in L^1([x_0,c))^{m\times m}$ for all $c>x_0$, and that $x$ is
a right Lebesgue point of $Q^{(N-1)}$. Let $M_+(z,x)$, $x\geq x_0$,
be either the unique limit point or a point of the limit disk $D_+(z,x)$
associated with \eqref{2.2}. Then, as $\abs{z}\to\infty$ in
$C_\varepsilon$, $M_+(z,x)$ has an asymptotic expansion of the form
$(\Im(z^{1/2})>0$, $z\in\bbC_+)$
\begin{equation}
M_+(z,x)\underset{\substack{\abs{z}\to\infty\\ z\in
C_\varepsilon}}{=}  i I_m z^{1/2}+\sum_{k=1}^N m_{+,k}(x)z^{-k/2}+
o(|z|^{-N/2}), \quad N\in\bbN. \lb{4.37}
\end{equation}
The expansion \eqref{4.37} is uniform with respect to $\arg\,(z)$ for 
$|z|\to \infty$ in $C_\varepsilon$ and uniform in $x$ as long as $x$
varies in compact subsets of $[x_0,\infty)$ intersected with the right
Lebesgue set of $Q^{(N-1)}$. The expansion coefficients $m_{+,k}(x)$ can
be recursively computed from \eqref{4.2} {\rm (}replacing $\tilde
m_{+,k}(x)$ by $m_{+,k}(x)${\rm )}. 

\noindent If one merely assumes Hypothesis~\ref{h2.1}, then \eqref{4.37} 
holds with $N=0$ {\rm (}interpreting $\sum_{k=1}^0 \cdot =0${\rm )},
uniformly with respect to $\arg (z)$ for $|z|\to \infty$ in
$C_\varepsilon$ and uniformly in $x$ as long as $x$ varies in compact
subsets of $[x_0,\infty)$. 
\end{theorem}

\begin{remark} \lb{r4.8}
Our asymptotic results in Theorem~\ref{t4.7} are not
necessarily confined to $M_+(z,x_0)$-matrices associated
with half-line Schr\"odinger operators $H_+$ on $[x_0,\infty)$.
In fact, introducing in addition the analogous Weyl-Titchmarsh
matrix $M_-(z,x_0)$ associated with a half-line Schr\"odinger
operator $H_-$ on $(-\infty,x_0]$, and noticing that the
diagonal Green's matrix $G(z,x,x)=(H-z)^{-1}(x,x)$ of a
matrix-valued Schr\"odinger  operator $H=-I_m\f{d^2}{dx^2}+Q$
in $L^2(\bbR)^m$ is given by
\begin{equation}
G(z,x,x)=(M_-(z,x)-M_+(z,x))^{-1}, \lb{4.38}
\end{equation}
Theorem~\ref{t4.6} then yields an analogous
asymptotic expansion for $G(z,x,x)$ in $C_\varepsilon$ of the form
\begin{equation}
G(z,x,x)\underset{\substack{\abs{z}\to\infty\\ z\in
C_\varepsilon}}{=}
\begin{cases}
 (i/2) I_m z^{-1/2}+o(|z|^{-1}) & \text{for $N=0$}, \\
 (i/2) \sum_{k=0}^N G_k(x) z^{-k-1/2}+o(|z|^{-N-1/2}) &
\text{for $N\in\bbN$}, \lb{4.39}
\end{cases}
\end{equation}
where
\begin{equation}
G_0(x)=I_m, \quad G_1(x)=\f12 Q(x), \text{ etc.} \lb{4.40}
\end{equation}
The expansion \eqref{4.39} is uniform with respect to
$\arg\,(z)$ for $|z|\to \infty$ in
$C_\varepsilon$ and uniform in $x\in\bbR$ as long as
$x$ varies in compact intervals interesected with the 
Lebesgue set of $Q^{(N-1)}$ (if $N\in\bbN$). We refer to
\cite{CGHL99} for further details and applications of this fact
in connection with trace formulas and Borg-type uniqueness
results for $Q(x)$, $x\in\bbR$.
\end{remark}

\vspace*{3mm}
\noindent {\bf Acknowledgments.}
We would like to thank Don Hinton, Helge Holden, Boris
Levitan, Mark Malamud, Alexei Rybkin, and Barry Simon for helpful
discussions and hints regarding the literature.




\begin{thebibliography}{99}

%
\bi{AM63} Z.~S.~Agranovich and V.~A.~Marchenko, {\it The Inverse
Problem of
Scattering Theory}, Gordon and Breach, New York, 1963.
%
\bi{AK92} M.~S.~Almamevov and V.~I.~Kantanova, {\it A regularized
trace of a high-order differential operator with a bounded
operational coefficient}, Diff. Eqs. {\bf 28}, 1--15 (1992).
%
\bi{AG98} D.~Alpay and I.~Gohberg, {\it Inverse problem for
Sturm-Liouville operators with rational reflection coefficient},
Integr. Equ. Oper. Theory {\bf 30}, 317--325 (1998).
%
\bi{AD56} N.~Aronszajn and W.~F.~Donoghue, {\it On
exponential representations of analytic functions in the
upper half-plane with positive imaginary part}, J. Anal.
Math. {\bf 5}, 321-388 (1956-57).
%
\bi{At64} F.~V.~Atkinson, {\it Discrete and Continuous Boundary
Problems},
Academic Press, New York, 1964.
%
\bi{At81} F.~V.~Atkinson, {\it On the location of Weyl circles},
Proc.
Roy. Soc. Edinburgh {\bf 88A}, 345--356 (1981).
%
\bi{At88} F.~V.~Atkinson, {\it Asymptotics of the Titchmarsh-Weyl
function in the matrix case}, unpublished preprint.
%
\bi{Be88} C.~Bennewitz, {\it A note on the Titchmarsh-Weyl
$m$-function},
Argonne Nat. Lab. preprint, ANL-87-26, Vol. 2, 1988, 105--111.
%
\bi{Be89} C.~Bennewitz, {\it Spectral asymptotics for
Sturm-Liouville equations}, Proc. London Math. Soc. (3) {\bf 59},
294--338 (1989).
%
\bi{Bo46} G.~Borg, {\it Eine Umkehrung der Sturm-Liouvilleschen
Eigenwertaufgabe}, Acta Math. {\bf 78}, 1--96 (1946).
%
\bi{BM97} A.~Boutet de Monvel and V.~Marchenko, {\it Asymptotic
formulas for spectral and Weyl functions of Sturm-Liouville
operators with smooth potentials}, in {\it New Results in Operator
Theory and its Applications}, I.~Gohberg and Yu.~Lubich (eds.),
Operator Theory, Advances and Applications, Vol.~98, Birkh\"auser,
Basel, 1997, pp.~102--117.
%
\bi{Ca76} R. W. Carey, {\it A unitary invariant for pairs
of self-adjoint
operators}, J. reine angew. Math. {\bf 283}, 294--312 (1976).
%
\bi{Ca98} R.~Carlson, {\it Large eigenvalues and trace formulas
for matrix Sturm-Liouville problems}, SIAM J. Math. {\bf 30}, 
949--962 (1999).
%
\bi{Ca98a} R.~Carlson, {\it Compactness of Floquet isospectral
sets for the matrix Hill's equation}, Proc. Amer. Math. Soc., 
{\bf 128}, 2933--2941 (2000).
%
\bi{Ca99} R.~Carlson, {\it Eigenvalue estimates and trace
formulas for the matrix Hill's equation}, J. Diff. Eq. {\bf 167}, 
211--244 (2000). 
%
\bi{Ch96} I.~Cherednik, {\it Basic Methods of Soliton
Theory}, World
Scientific, Singapore, 1996.
%
\bi{Ch99} H.-H.~Chern, {\it On the eigenvalues of some
vectorial Sturm-Liouville eigenvalue problems}, preprint, 1999.
%
\bi{Ch99a} H.-H.~Chern, {\it On the construction of isospectral
vectorial Sturm-Liouville differential equations}, preprint, 1998.
%
\bi{CS97} H.-H.~Chern and C-L.~Shen, {\it On the $n$-dimensional
Ambarzumyan's theorem}, Inverse Problems {\bf 13}, 15--18 (1997).
%
\bi{Cl92} S.~L.~Clark, {\it Asymptotic behavior of the
Titchmarsh-Weyl coefficient for a coupled second order system}, 
in ``{\it Ordinary and Delay Differential Equations\/}'', J.~Wiener 
and J.~K.~Hale (eds.), Longman, New York, 1992, pp.~24--28.
%
\bi{CG00} S.~Clark and F.~Gesztesy, {\it Weyl-Titchmarsh
$M$-function asymptotics and Borg-type theorems for Dirac-type
operators}, in preparation.
%
\bi{CGHL99} S.~Clark, F.~Gesztesy, H.~Holden, and B.~M.~Levitan,
{\it Borg-type theorems for matrix-valued Schr\"odinger
operators}, J. Diff. Eq. {\bf 167}, 181--210 (2000).
%
\bi{CGR99} S.~Clark, F.~Gesztesy, and W.~Renger, Trace formulas and 
Borg-type theorems for finite difference operators, in preparation.
%
\bi{Co77} M.~Coz, {\it The Riemann solution in the one-dimensional
inverse problem}, J. Math. Anal. Appl. {\bf 61}, 232--250 (1977).
%
\bi{DL91} A.~A.~Danielyan and B.~M.~Levitan, {\it On the asymptotic
behavior of the Weyl-Titchmarsh $m$-function}, Math. USSR Izv.
{\bf 36}, 487--496 (1991).
%
\bi{De95} B.~Despr\'es, {\it The Borg theorem for the vectorial
Hill's equation}, Inverse Probl. {\bf 11}, 97--121 (1995).
%
\bibitem{Di91}  L.~A.~Dickey, {\it Soliton Equations and
Hamiltonian
Systems},  World Scientific, Singapore, 1991.
%
\bi{DK97} V.~S.~Dryuma and B.~G.~Konopelchenko, {\it On equation
of geodesic deviation and its solution}, Izv. Akad. Nauk Respub. 
Moldova Mat. no. 3, 61--73, 121, 123, (1996).
%
\bibitem{Du77}
B.~A.~Dubrovin, {\it Completely integrable Hamiltonian Systems
associated with Matrix operators and Abelian varieties},
Funct. Anal. Appl. \textbf{11}, 265--277 (1977).
%
\bibitem{Du83}
B.~A.~Dubrovin, {\it  Matrix finite-zone operators},
Revs. Sci. Tech. \textbf{23}, 20--50 (1983).
%
\bi{Ev72} W.~N.~Everitt, {\it On a property of the $m$-coefficient
of a second-order linear differential equation}, J. London Math. 
Soc. (2), {\bf 4}, 443--457 (1972).
%
\bi{EH78} W.~N.~Everitt and S.~G.~Halvorsen, {\it On the asymptotic
form of the Titchmarsh-Weyl $m$-coefficient}, Appl. Anal. {\bf 8},
153--169 (1978).
%
\bi{EHS83} W.~N.~Everitt, D.~B.~Hinton, and J.~K.~Shaw, {\it The
asymptotic form of the Titchmarsh-Weyl coefficient for Dirac 
systems}, J. London Math. Soc. (2), {\bf 27}, 465--476 (1983).
%
\bi{GD77} I.~M.~Gel'fand and L.~A.~Dikii, {\it The resolvent and
Hamiltonian systems}, Funct. Anal. Appl. {\bf 11}, 93--105 (1977).
%
\bi{GH97} F.~Gesztesy and H.~Holden, {\it On trace formulas for
Schr\"odinger-type operators}, in ``{\it Multiparticle Quantum
Scattering with Applications to Nuclear, Atomic and Molecular
Physics\/}'', D.~G.~Truhlar and B.~Simon (eds.), Springer, New 
York, 1997, p.~121--145.
%
\bi{GHSZ95} F.~Gesztesy, H.~Holden, B.~Simon, and Z.~Zhao,
{\it Higher order trace relations for Schr\"odinger
operators}, Rev. Math. Phys. {\bf 7}, 893--922 (1995).
%
\bi{GKMT98} F.~Gesztesy, N.~J.~Kalton, K.~A.~Makarov,
and E.~Tsekanovskii,
{\it Some Applications of Operator-Valued Herglotz
Functions}, Operator Theory: Advances and Applications,
Birkh\"auser (to appear).
%
\bi{GM99} F.~Gesztesy and K.~A.~Makarov, {\it Some
applications of the spectral shift operator}, in
{\it Operator Theory and its Applications}, A.~G.~Ramm,
P.~N.~Shivakumar, and A.~V.Strauss (eds.), Fields
Institute Communications, Vol.~25, Amer. Math. Society,
Providence RI, 2000, p.~267--292.
%
\bi{GM99a} F.~Gesztesy and K.~A.~Makarov, {\it The
$\Xi$ operator and its relation to Krein's spectral
shift function}, J. d'Anal. Math. {\bf 81}, 139--183 (2000).
%
\bi{GMN98} F.~Gesztesy, K.~A.~Makarov, and S.~N.~Naboko, {\it The
spectral shift operator}, in {\it Mathematical Results in Quantum
Mechanics}, J.~Dittrich, P.~Exner, and M.~Tater (eds.), Operator
Theory: Advances and Applications,  Vol. 108, Birkh\"auser, Basel, 
1999, p.~59--90.
%
\bi{GMT98} F.~Gesztesy, K.~A.~Makarov, and E.~Tsekanovskii,
{\it An Addendum to Krein's Formula,}
J. Math. Anal. Appl. {\bf 222}, 594--606 (1998).
%
\bi{GS96} F.~Gesztesy and B.~Simon, {\it The xi function,}
Acta Math. {\bf 176}, 49--71 (1996).
%
\bi{GS98} F.~Gesztesy and B.~Simon, {\it A new approach to inverse
spectral theory, II. General real potentials and the connection to
the spectral measure}, Ann. of Math. {\bf 152}, 593--643 (2000).
%
\bi{GT97} F.~Gesztesy and E.~Tsekanovskii,
{\it On matrix-valued Herglotz functions,} Math. Nachr. {\bf 218}, 
61--138 (2000).
%
\bi{GKS97} I.~Gohberg, M.~A.~Kaashoek, and A.~L.~Sakhnovich, {\it
Sturm-Liouville systems with rational Weyl functions: explicit
formulas and applications}, Integral Eq. Operator Th. {\bf 30},
338--377, (1998).
%
\bi{Ha83} B.~J.~Harris, {\it On the Titchmarsh-Weyl $m$-function},
Proc. Roy. Soc. Edinburgh {\bf 95A}, 223--237 (1983).
%
\bi{Ha84} B.~J.~Harris, {\it The asymptotic form of the
Titchmarsh-Weyl
$m$-function}, J. London Math.
Soc. (2), {\bf 30}, 110--118 (1984).
%
\bi{Ha86} B.~J.~Harris, {\it The asymptotic form of the
Titchmarsh-Weyl
$m$-function associated with a second order differential equation
with
locally integrable coefficient}, Proc. Roy. Soc. Edinburgh
{\bf 102A},
243--251 (1986).
%
\bi{Ha87} B.~J.~Harris, {\it An exact method for the calculation
of certain Titchmarsh-Weyl $m$-functions},  Proc. Roy. Soc.
Edinburgh {\bf 106A},
137--142 (1987).
%
\bi{Hi69} E.~Hille, {\it Lectures on Ordinary Differential
Equations}, Addison-Wesley, Reading, 1969.
%
\bi{HKS89} D.~B.~Hinton, M.~Klaus, and J.~K.~Shaw, {\it Series
representation and asymptotics for Titchmarsh-Weyl $m$-functions},
Diff. Integral Eqs. {\bf 2}, 419--429 (1989).
%
\bi{HSH93}  D.~B.~Hinton and A. Schneider, {\it On the Titchmarsh-Weyl
coefficients for singular S-Hermitian Systems I}, Math. 
Nachr. {\bf 163}, 323--342 (1993).
%
\bi{HSH97}  D.~B.~Hinton and A. Schneider, {\it On the Titchmarsh-Weyl
coefficients for singular S-Hermitian Systems II}, Math. 
Nachr. {\bf 185}, 67--84 (1997).
%
\bi{HS81} D.~B.~Hinton and J.~K.~Shaw, {\it On Titchmarsh-Weyl
$M(\lambda)$-functions for linear Hamiltonian systems}, J. Diff.
Eqs. {\bf 40}, 316--342 (1981).
%
\bi{HS82} D.~B.~Hinton and J.~K.~Shaw, {\it On the spectrum of a
singular Hamiltonian system}, Quaest. Math. {\bf 5}, 29--81 (1982).
%
\bi{HS83} D.~B.~Hinton and J.~K.~Shaw, {\it Hamiltonian systems of
limit point or limit circle type with both endpoints singular}, J.
Diff. Eqs. {\bf 50}, 444--464 (1983).
%
\bi{HS84} D.~B.~Hinton and J.~K.~Shaw, {\it On boundary value
problems for Hamiltonian systems with two singular points}, SIAM
J. Math. Anal. {\bf 15}, 272--286 (1984).
%
\bi{HS86} D.~B.~Hinton and J.~K.~Shaw, {\it On the spectrum of a
singular Hamiltonian system, II}, Quaest. Math. {\bf 10}, 1--48
(1986).
%
\bi{JL97} M.~Jodeit and B.~M.~Levitan, {\it Isospectral
vector-valued Sturm-Liouville problems}, Lett. Math. Phys.
{\bf 43}, 117--122 (1998).
%
\bi{Jo87} R.~A.~Johnson, {\it $m$-Functions and Floquet exponents
for linear
differential systems}, Ann. Mat. Pura Appl., Ser. 4, {\bf 147},
211--248
(1987).
%
\bi{KK86} H.~G.~Kaper and M.~M.~Kwong, {\it Asymptotics of the
Titchmarsh-Weyl $m$-coefficient for integrable potentials}, Proc.
Roy. Soc. Edinburgh {\bf 103A}, 347--358 (1986).
%
%
\bi{KK87} H.~G.~Kaper and M.~M.~Kwong, {\it Asymptotics of the
Titchmarsh-Weyl $m$-coefficient for integrable potentials, II}, in
{\it Differential Equations and Mathematical Physics},
I.~W.~Knowles and Y.~Saito (eds.), Lecture Notes in Mathematics,
Vol.~1285, Springer, Berlin, 1987, pp.~222--229.
%
\bi{KR74} V.~I.~Kogan and F.~S.~Rofe-Beketov, {\it On
square-integrable
solutions of symmetric systems of differential equations of arbitrary
order}, Proc. Roy. Soc. Edinburgh {\bf 74A}, 1--40 (1974).
%
\bi{KS88} S.~Kotani and B.~Simon,
{\it Stochastic Schr\"odinger operators and Jacobi matrices on
the strip}, Commun. Math. Phys. {\bf 119},  403--429 (1988).
%
\bi{Kr89a} A.~M.~Krall, {\it $M(\lambda)$ theory for singular
Hamiltonian
systems with one singular point}, SIAM J. Math. Anal. {\bf 20},
664--700 (1989).
%
\bi{Kr89b} A.~M.~Krall, {\it $M(\lambda)$ theory for singular
Hamiltonian
systems with two singular points}, SIAM J. Math. Anal. {\bf 20},
701--715 (1989).
%
\bi{Kr83} M.~G.~Krein, {\it On tests for stable boundedness of
solutions of periodic canonical systems}, Amer. Math. Soc. Transl.
(2) {\bf 120}, 71--110 (1983).
%
\bi{Ma94} M.~M.~Malamud, {\it Similarity of Volterra operators and
related questions of the theory of differential equations of
fractional
order}, Trans. Moscow Math. Soc. {\bf 55}, 57--122 (1994).
%
\bi{Ma98} M.M.~Malamud,
{\it Uniqueness questions in inverse problems for systems of
ordinary differential equations on a finite interval}, Trans. Moscow
Math. Soc.  {\bf 60}, 173--224  (1999).
%
\bi{Ma78} Yu.~I.~Manin, {\it Matrix solitons and bundles over
curves with
singularities}, Funct. Anal. Appl. {\bf 12}, 286--295 (1978).
%
\bi{Ma88} V.~A.~Marchenko, {\it Nonlinear Equations and Operator
Algebras}, Reidel, Dordrecht, 1988.
%
\bi{MO82} L.~Martinez Alonso and E.~Olmedilla, {\it Trace
identities in the inverse
scattering transform method associated with matrix Schr\"odinger
operators}, J.
Math. Phys. {\bf 23}, 2116--2121 (1982).
%
\bi{MPS90} K.~N.~Murty, K.~R.~Prasad, and M.~A.~S.~Srinivas, {\it
Upper and lower bounds for the solution of the general matrix
Riccati differential equation}, J. Math. Anal. Appl. {\bf 147},
12--21 (1990).
%
\bi{NJ55} R.~G.~Newton and R.~Jost, {\it The construction of
potentials from the $S$-matrix for systems of differential
equations}, Nuovo Cim. {\bf 1}, 590--622 (1955).
%
\bi{OMG81} E.~Olmedilla, L. Martinez Alonso, and F. Guil, {\it
Infinite-dimensional Hamiltonian systems associated with matrix
Schr\"odinger operators}, Nuovo Cim. {\bf 61 B}, 49--61 (1981).
%
\bi{Or76} S.~A.~Orlov, {\it Nested matrix disks analytically
depending
on a parameter, and theorems on the invariance of ranks of radii of
limiting disks}, Math. USSR Izv. {\bf 10}, 565--613 (1976).
%
\bi{Pa95} V.~G.~Papanicolaou, {\it Trace formulas and the behavior
of large eigenvalues}, SIAM J. Math. Anal. {\bf 26}, 218--237
(1995).
%
\bi{Re96} C.~Remling, {\it Geometric characterization of singular
self-adjoint boundary conditions for Hamiltonian systems}, Appl.
Anal. {\bf 60}, 49--61 (1996).
%
\bi{Ro60} F.~S.~Rofe-Beketov, {\it Expansions in eigenfunctions of
infinite systems of differential equations in the non-self-adjoint
and self-adjoint case}, Mat. Sb. {\bf 51}, 293--342 (1960)
(Russian).
%
\bi{Ro69} F.~S.~Rofe-Beketov, {\it Selfadjoint extensions of
differential operators in a space of vector functions}, Sov. Math.
Dokl. {\bf 10}, 188--192 (1969).
%
\bi{Ry99} A.~Rybkin, {\it On the trace approach to 
the inverse scattering problem in dimension one}, SIAM J. Math. 
Anal. (to appear).
%
\bi{Sa94} L.~A.~Sakhnovich, {\it Inverse problems for equations
systems}, in {\it Matrix and Operator Valued Functions: The Vladimir
Petrovich Potapov Memorial Volume}, I.~Gohberg and L.~A.~Sakhnovich
(eds.),  Operator Theory: Advances and Applications, Vol. 72,
Birkh\"auser, Basel, 1994, p.~202--211.
%
\bi{Sa94a}  L.~A.~Sakhnovich, {\it Method of operator identities
and problems of analysis}, St. Petersburg Math. J. {\bf 5}, 1--69
(1994).
%
\bi{Sc83} P.~Schuur, {\it Inverse scattering for the matrix
Schr\"odinger
equation with non-hermitian potential}, in ``{\it Nonlinear
Waves\/}'',
L.~Debnath (ed.), Cambridge University Press, Cambridge, 1983,
pp.~285--297.
%
\bi{Sh71} Yu.~L.~Shmul'yan,
{\it On operator $R$-functions,}
Siberian Math. J. {\bf 12}, 315--322  (1971).
%
\bi{Si98} B.~Simon, {\it A new approach to inverse spectral theory,
I. Fundamental formalism}, Ann. of Math. {\bf 150}, 1029--1057 
(1999).
%
\bi{Ti86} E.~C.~Titchmarsh, {\it Introduction to the Theory of Fourier
Integrals}, Chelsea, New York, 1986.
%
\bi{WK74} M.~Wadati and T.~Kamijo, {\it On the extension of inverse
scattering method}, Progr. Theoret. Phys. {\bf 52}, 397--414 (1974).
%
\bi{We87} J.~Weidmann, {\it Spectral Theory of Ordinary Differential
Operators}, Lecture Notes in Math. {\bf 1258}, Springer, Berlin,
1987.
%

\end{thebibliography}
\end{document}